\documentclass[12pt,draftcls,onecolumn]{IEEEtran}
\IEEEoverridecommandlockouts

\usepackage{fullpage}
\usepackage{amsmath}
\usepackage{graphicx}
\usepackage{amssymb}
\usepackage{color}
\usepackage{authblk}
\usepackage{psfrag}
\begin{document}
\title{\LARGE{Periodic Event-Triggered Synchronization of Linear Multi-agent Systems with Communication Delays}}
\author{Eloy Garcia\thanks{Corresponding Author: \ttfamily{elgarcia@infoscitex.com}. \\ This work has been supported in part by AFOSR LRIR No. 12RB07COR. \\ Eloy Garcia is a contractor (Infoscitex Corp.) with the Control Science Center of Excellence, Air Force Research Laboratory, Wright-Patterson AFB, OH 45433. \\ Yongcan Cao and David Casbeer are with the Control Science Center of Excellence, Air Force Research Laboratory, Wright-Patterson AFB, OH 45433. \\ A shorter version of this document was submitted to the IEEE Transactions on Automatic Control.}, Yongcan Cao, and David W. Casbeer}          
         
\newtheorem{theorem}{Theorem}
\newtheorem{lemma}{Lemma}

\maketitle 
\begin{abstract}
Multi-agent systems cooperation to achieve global goals is usually limited by sensing, actuation, and communication issues. At the local level, continuous measurement and actuation is only approximated by the use of digital mechanisms that measure and process information in order to compute and update new control input values at discrete time instants. Interaction with other agents or subsystems takes place, in general, through a digital communication channel with limited bandwidth where transmission of continuous-time signals is not possible. Additionally, communication channels may be subject to other imperfections such as time-varying delays. This paper considers the problem of consensus (or synchronization of state trajectories) of multi-agent systems that are described by general linear dynamics and are connected using undirected graphs. An event-triggered consensus protocol is proposed, where each agent implements discretized and decoupled models of the states of its neighbors. This approach not only avoids the need for continuous communication between agents but also provides a decentralized method for transmission of information in the presence of time-varying communication delays where each agent decides its own broadcasting time instants based only on local information. This method gives more flexibility for scheduling information broadcasting compared to periodic and sampled-data implementations. The use of discretized models by each agent allows for a periodic event-triggered strategy where continuous actuation and continuous measurement of the states are not necessary.  
\end{abstract}

\section{Introduction} \label{sec:one}

Cooperative control of multi-agent systems is an active research area with broad and relevant applications in commercial, academic and military areas~\cite{RenBeard07}. The design of decentralized and scalable control algorithms provides the necessary coordination for a group of agents to outperform agents operating independently. In multiple scenarios where communication among agents is limited, decentralized computation of the time instants that each agent needs to transmit relevant information is also necessary. Continuous actuation and continuous measurement of local states may also be restricted by particular hardware limitations.


Consensus problems with limited communication have been studied using the sampled-data (periodic) approach \cite{CaoRen10}, \cite{Hayakawa06}, \cite{Liu10}, and \cite{QinGao12}. An important drawback of periodic transmission is that it requires synchronization between the agents, that is, all agents need to transmit their information at the same time instants and, in some cases, it requires a conservative sampling period for worst case situations. 

In the present paper, in lieu of periodic approaches, we use an asynchronous communication scheme based on event-triggered control strategies and we consider agents that are described by general linear dynamics which are subject to limited actuation update rates and also to limited local sensor measurement update rates. In addition, we consider the case where communication among agents is subject to communication delays. In contrast to periodic (or time-triggered) implementations, in the context of event-triggered control, information or measurements are not transmitted periodically in time but they are triggered by the occurrence of certain events. In event-triggered broadcasting  \cite{Astrom08}, \cite{Donkers10}, \cite{DePersis13}, \cite{Garcia13}, \cite{Wang11}, \cite{Garcia12IJC}, \cite{Garcia12CDC}, and \cite{Tabuada07}, a subsystem sends its local state to the network only when it is necessary, that is, only when a measure of the local subsystem state error is above a specified threshold. Event-triggered control strategies have been used for stabilization of multiple coupled subsystems as in \cite{Garcia12}, \cite{Stocker13}, and \cite{Guinaldo11}. Consensus problems have also been studied using these techniques \cite{Dimarogonas12}, \cite{Garcia13b}, \cite{Seyboth13}, \cite{Yu12}, \cite{YinYue13b}, \cite{ChenHao12}, \cite{GuoDimarogonas13}. Event-triggered control provides a more robust and efficient use of network bandwidth. Its implementation in multi-agent systems also provides a highly decentralized way to schedule transmission instants which does not require synchronization compared to periodic sampled-data approaches. 

One important restriction related to event-triggered control techniques is that continuous measurement of state variables and continuous computation of state errors and time-varying thresholds is required. One solution explored by different researchers is self-triggered control \cite{AntaTabuada10}, \cite{Mazo10}, \cite{Nowzari12}, \cite{WangLemmon09}. The main difference with respect to event-triggered control is that a measure of the state is not being compared constantly against a predefined threshold. Instead, the current state measurement is used to determine its next deadline, i.e. the next time that the sensor is required to send a measurement to the controller. A recent extension to event-triggered control is the so called \textit{periodic event-triggered control} \cite{Heemels13periodic}, \cite{Heemels13model}, where measurements of states and computations of errors and thresholds occur, not continuously in time, but only at periodic time instants. The decision to transmit the current measurement is taken at the sampling instants based on current measurements and computed thresholds.

Consensus problems where all agents are described by general linear models have been considered by different authors \cite{LiDuan10}, \cite{LiDuan11}, \cite{MaZhang10}, \cite{Ren08}, \cite{Scardovi09}, \cite{SuHuang12}, \cite{Tuna08}, and \cite{Tuna09}. In these papers it is assumed that continuous communication between agents is possible. The work in \cite{WenDuan13} considers the consensus problem of agents with linear dynamics under communication constraints. Specifically, the authors consider the existence of continuous communication among agents for finite intervals of time and the total absence of communication among agents for other time intervals, and the minimum rate of continuous communication to no communication is given.

Event-triggered consensus of agents with linear dynamics and limited communication was recently explored in \cite{LiuHill12} and \cite{Zhu14}. In our previous work \cite{Garcia14} and \cite{Garcia14auto} we proposed a novel approach in which each agent implements models of the decoupled dynamics of each one of its neighbors and uses the model states to compute the local control input. This approach offered better performance than Zero-Order-Hold (ZOH) approach used in \cite{LiuHill12} and \cite{Zhu14} where the updates from neighbors are kept constant by the local agent. A similar model-based framework was proposed in \cite{Demir12b} where only constant thresholds were used. One of the main limitations of the ZOH approach \cite{LiuHill12}, \cite{Zhu14} is that it is not capable to keep up with unstable trajectories and updates need to be generated more frequently. In consensus with general linear dynamics, unstable systems are one of the most interesting cases to analyze. The model-based approach in \cite{Garcia14}, \cite{Garcia14auto} provides better estimates of neighbors and reduces generation of events as agents converge to similar unstable trajectories. Communication delays using the event-triggered approach in these papers were addressed in \cite{Garcia14CDC}. The present paper offers complete proofs with respect to \cite{Garcia14CDC} and, more importantly, it extends the consensus protocol in that paper in order to consider limited actuation and sensing update rates. A periodic event-triggered control technique for consensus of linear multi-agent systems is proposed in the present paper where each agent implements discretized and decoupled models of the states of neighbors. Communication delays are also considered. The case shown in \cite{Garcia14CDC} represents only a particular case of the approach described here, when continuous actuation and continuous measurements are possible to implement and obtain, respectively.

In the present paper decentralized event thresholds that guarantee practical consensus\footnote{In the case of unstable trajectories practical consensus is obtained when the difference between the states of any two agents is bounded.} and strictly positive inter-event times are designed. The lower-bounds on the inter-event time intervals are independent of the particular system trajectories, therefore they hold for any two consecutive local events. The main contribution of this paper is the design of periodic event-triggered strategies for consensus of linear systems with limited sensing and actuation updates and with communication delays. The periodic event-triggered control technique automatically avoids the presence of Zeno behavior\footnote{Zeno behavior in event-triggered control refers to the occurrence of an infinite number of triggering events in a finite time interval.}. However, and for completeness, we establish the relationship between a selected sampling period and the performance of the consensus protocol with respect to the bounds on the state disagreement.

The present paper addresses the several problems that were not considered by \cite{Garcia14auto}:
\begin{enumerate}
	\item Time-varying communication delays. The transmitted event-based measurement updates are subject to communication delays. 
  \item Constrained sensing and event computation rate. Each agent does not need to continuously measure its own state but only at finite time instants. Similarly, continuous computation of errors and events is not necessary and these operations are performed only at discrete-time instants.
  \item Constrained actuation rate. This event-based approach also provides sampled actuation time instants instead of continuous actuation.
\end{enumerate}
Continuous measurement, actuation, and computation of events severely restrict the operation of the subsystems; therefore, the relaxations addressed in this paper offer a significant advantage in terms of implementation and resource management. 

The remainder of this paper is organized as follows. Section~\ref{sec:two} provides a brief background on graph theory and describes the problem and the consensus protocol. Section~\ref{sec:three} gives a result assuming continuous communication which will be used in the main results of this paper in Sections~\ref{sec:four} and~\ref{sec:five}. Design of periodic decentralized event thresholds for systems with limited sensing and actuation capabilities is addressed in Section~\ref{sec:four}. Section~\ref{sec:five} extends this approach in order to consider time-varying, but bounded, communication delays. Section~\ref{sec:six} presents illustrative examples and Section~\ref{sec:seven} concludes the paper.

\section{Preliminaries} \label{sec:two}
\subsection{Graph Theory}
Consider a graph $\mathcal{G}=\left\{\mathcal{V,E}\right\}$ consisting of a set of vertices or nodes $\mathcal{V}=\left\{1,...,N\right\}$ and a set of edges $\mathcal{E}$. An edge between nodes $i$ and $j$ is represented by the pair $(i,j)\in \mathcal{E}$. A graph $\mathcal{G}$ is called undirected if $(i,j)\in \mathcal{E}\Leftrightarrow (j,i)\in \mathcal{E}$ and the nodes are called adjacent. The adjacency matrix $\mathcal{A}$ is defined by $a_{ij}=1$ if the nodes $i$ and $j$ are adjacent and $a_{ij}=0$ otherwise. If $(j,i)\in \mathcal{E}$, then $j$ is said to be a neighbor of $i$. The set $\mathcal{N}_i$ is called the set of neighbors of node $i$, and $N_i$ is its cardinality. A node $j$ is an element of $\mathcal{N}_i$ if $(j,i)\in \mathcal{E}$. 
A path from node $i$ to node $j$ is a sequence of distinct nodes that starts at $i$ and ends at $j$, such that every pair of consecutive nodes is adjacent. An undirected graph is connected if there is a path between every pair of distinct nodes. The Laplacian matrix        
$\mathcal{L}$ of $\mathcal{G}$ is defined as $\mathcal{L=D-A}$ where $\mathcal{D}$ represents the degree matrix which is a diagonal matrix with entries $d_{ii}=\sum_{j\in\mathcal{N}_i}a_{ij}$. For undirected graphs, $\mathcal{L}$ is symmetric and positive semi-definite. $\mathcal{L}$ has zero row sums and, therefore, zero is an eigenvalue of $\mathcal{L}$ with associated eigenvector $\textbf{1}_N$ (a vector with all its $N$ entries equal to one), that is, $\mathcal{L}\textbf{1}_N=\textbf{0}_N$. If an undirected graph is connected then $\mathcal{L}$ has exactly one eigenvalue equal to zero and all its non-zero eigenvalues are positive; they can be set in increasing order $\lambda_1(\mathcal{L})<\lambda_2(\mathcal{L})\leq \lambda_3(\mathcal{L})\leq ...\leq \lambda_N(\mathcal{L})$, with $\lambda_1(\mathcal{L})=0$. 

\begin{lemma} \label{lm:LyapFn}
Let $\mathcal{L}$ be the symmetric Laplacian of an undirected and connected graph. Then, consensus is achieved if and only if 
\begin{align}
	V=\chi^T\hat{\mathcal{L}}\chi =0, \label{eq:LyapV}
\end{align} 
where $\hat{\mathcal{L}}=\mathcal{L}\otimes Q$, $Q\in \mathbb{R}^{n\times n}$ is a symmetric positive definite matrix,
\begin{align}
\chi(t)=\begin{bmatrix}
	\chi_1(t)^T&\chi_2(t)^T&\ldots&\chi_N(t)^T  \nonumber
\end{bmatrix}^T,
\end{align}
$\chi_i\in \mathbb{R}^n$, and $\otimes$ denotes the Kronecker product.
\end{lemma}

\textit{Proof}. See Appendix.

\subsection{Problem Statement}
\begin{figure}
	\begin{center}
		\includegraphics[width=8.4cm,height=6.5cm,trim=1cm .1cm .8cm .1cm]{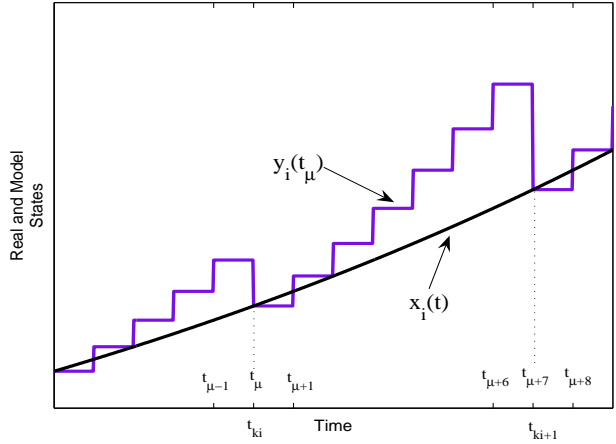}
	\caption{Representation of periodic sampling time instants and event-triggered communication time instants}
	\label{fig:tutki}
	\end{center}
\end{figure}

Consider a group of $N$ agents with fixed communication graphs and fixed weights. Each agent's dynamics are described by the following:
\begin{align}
	\dot{x}_i(t)=Ax_i(t)+Bu_i(t_\mu),\ \ i=1,...,N, \label{eq:agents}
\end{align}
 with
\begin{align}
	u_i(t_\mu)=cF\sum_{j\in\mathcal{N}_i}(x_i(t_\mu)-y_j(t_\mu)),\ \ i=1,...,N, \label{eq:inputs}
\end{align}
where $A\in \mathbb{R}^{n\times n}$, $B\in \mathbb{R}^{n\times m}$, $x_i\in \mathbb{R}^n$ is the state of agent $i$, and $u_i\in \mathbb{R}^m$ is the control input for agent $i$. $F\in\mathbb{R}^{m\times n}$ and $c\in\mathbb{R}_+$ are design parameters that are defined below. 
The variables $y_j\in \mathbb{R}^n$ represent a model of the $j^{th}$ agent's state using the decoupled and discretized dynamics:  
\begin{align}
	y_j(t_{\mu+1})=Gy_j(t_\mu),\ \ j=1,...,N. \label{eq:models}
\end{align}
for $t_\mu \in[t_{k_j},t_{k_j+1})$ where $G=e^{Ah}$, $h=t_{\mu+1}-t_\mu$, and $y_j(t_{k_j})=x_j(t_{k_j})$. We refer to $h$ as the discretization or sampling period since it is used to obtain the discrete-time model $G$ and the state of each agent, $x_i$, $i=1,...,N$, is sampled every $h$ time units. The discretization period $h$ is constant but the communication intervals for each agent are not constant. The transmission time intervals are non-periodic and are determined by event-triggered rules. This means that at every time instant $t_\mu$ each agent samples its own state and updates its control input \eqref{eq:inputs}. It also uses its sampled state to compute its local state error and to determine if it is necessary or not to transmit the current state $x(t_\mu)$ to its neighbors. An example showing this relationship is shown in Fig. \ref{fig:tutki}. This figure shows the continuous-time state $x_i(t)$ and its discrete-time model $y_i(t_\mu)$. Two event-triggered model updates are shown, at time $t_{k_i}$ and at time $t_{k_i+1}$. The first event-triggered model update occurs at sampling time $t_\mu$ while the following event time instant is equivalent, in this example, to sampling time $t_{\mu+7}$, that is, $7$ sampling time instants occur between $t_{k_i}$ and $t_{k_i+1}$ in this example. Note that at the event time instants the model is updated and it takes the current value of $x_i(t_{k_i})$.

The main advantage of using a periodic event-triggered control scheme is that measuring the state of each agent (and the associated computations that require evaluation of state errors and thresholds) is only performed at some periodic time instants (every $h$ time units) instead of doing it continuously as it is common in most event-triggered control schemes. There are additional advantages of using discrete-time models of the continuous-time system dynamics compared to \cite{Garcia14auto}, \cite{Garcia14CDC} where continuous-time models were used. For instance, continuous actuation is not required and the operations related to implementing discrete-time models at every node are simplified with respect to implementing continuous-time models.   

Every agent in the network implements a model of itself $y_i(t_\mu)$ and also models of its neighbors $y_j(t_\mu)$, $j\in\mathcal{N}_i$. The model state $y_i(t_\mu)$ is not used by agent $i$ for control since the real state, $x_i(t_\mu)$, is locally available at every sampling time $t_\mu$. However, the local model state, $y_i(t_\mu)$, is used to trigger local events. This way of defining the local control input \eqref{eq:inputs} represents a difference with respect to the approach in \cite{Garcia14} and \cite{Garcia14auto} where the control input is a function of model variables, including $y_i$.  

Local events for agent $i$ are defined as follows. When an event is triggered, agent $i$ will transmit its current state $x_i(t_{k_i})$ to its neighbors and each one of the neighbors updates the model state they have of agent $i$. At the same time instant, agent $i$ will update its local model using the same state measurement. The notation $t_{k_i}$ represents the local broadcasting time instants, i.e. when agent $i$ triggers an event and transmits its current state. Since we check for events not continuously but only at sampling times $t_\mu$, then $t_{k_i}$ is equal to some sampling instant $t_\mu$ but, in general, not every time instant $t_\mu$ is a broadcasting instant $t_{k_i}$. The decision to trigger an event in order to broadcast the current state is given by the event-triggered schemes described in the following sections.

For negligible communication delays, agent $i$ and its neighbors will all update their local models $y_i(t_\mu)$ at the same time instant. Since agent $i$ and its neighbors use the same measurements to update the models, say, $x_i(t_{k_i})$ and the model dynamics~\eqref{eq:models} represent the decoupled dynamics where all agents use the same state matrix, then the model states $y_i(t_\mu)$ implemented by agent $i$ and by its neighbors are the same. In the presence of communication delays the previous statement will not hold and we will differentiate between $y_{ii}(t_\mu)$, the model state of agent $i$ as seen by agent $i$, and $y_{ij}(t_\mu)$, the model state of agent $i$ as seen by agents $j$, $j\in \mathcal{N}_i$. More details concerning communication delays are presented in Section \ref{sec:five}.

The model update process is similar for all agents $i=1,...,N$. The local control input~\eqref{eq:inputs} is decentralized since it only depends on local information, that is, on the sampled state of the local agent and on the discretized model states of its neighbors. Continuous or even periodic access to the states of neighbors is not needed. 

Note that the difference between ZOH periodic samples of the agent dynamics~\eqref{eq:agents} and our proposed models~\eqref{eq:models} is given by the input term in~\eqref{eq:agents} and this input decreases as the agents approach a consensus state. It can also be seen that in the particular case when systems~\eqref{eq:agents} represent single integrator dynamics, then our models degenerate to ZOH models as in \cite{Dimarogonas12}, \cite{Garcia13b}.

\section{Consensus with continuous measurements} \label{sec:three}

This section provides an important result in Lemma \ref{lm:ConsLin} that will be useful in the following sections. Let us assume in this section that continuous actuation by each agent and continuous communication between agents are possible, then the agent dynamics \eqref{eq:agents} are now driven by the local control inputs
\begin{align}
	u_i(t)=cF\sum_{j\in\mathcal{N}_i}(x_i(t)-x_j(t)),\ \ i=1...N. \label{eq:ContInput}
\end{align}
Assume that the pair ($A,B$) is controllable. Then, for $\alpha>0$ there exists a symmetric and positive definite solution $P$ to
\begin{align}
	PA+A^TP-2PBB^TP+2\alpha P<0. \label{eq:LMI}
\end{align}
Let  
\begin{align}
  F&=-B^TP \label{eq:F}
\end{align}
and	$c\geq1/\lambda_2$.

\begin{theorem}\label{th:consensus}
 Assume the pair ($A,B$) is controllable and the communication graph is connected and undirected. Define $F$ as in~\eqref{eq:F} and	$c\geq1/\lambda_2$. Then the following symmetric matrix
\begin{align}
  \bar{\mathcal{L}}=\hat{\mathcal{L}}A_c+A^{T}_c\hat{\mathcal{L}} \label{eq:Lbar}
\end{align}
has only $n$ eigenvalues equal to zero and the rest of its eigenvalues are negative. In addition, the eigenvectors associated with its $n$ zero eigenvalues belong to the subspace spanned by the eigenvectors associated with the $n$ zero eigenvalues of $\hat{\mathcal{L}}$, 
where $\hat{\mathcal{L}}=\mathcal{L}\otimes P$, $A_c=\bar{A}+\bar{B}$, $\bar{A}=I_N\otimes A$, $\bar{B}=c\mathcal{L}\otimes BF$, and $I_N$ is an identity matrix of size $N$.
\end{theorem}

\textit{Proof}. See Appendix. 

\begin{lemma} \label{lm:ConsLin}
Assume the pair ($A,B$) is controllable and the communication graph is connected and undirected. Then, protocol \eqref{eq:ContInput}, with $F$ defined in~\eqref{eq:F} and $c\geq1/\lambda_2$, solves the consensus problem for agents described by \eqref{eq:agents}. Furthermore, the Lyapunov function defined by $V=x^{T}\hat{\mathcal{L}}x$ has a time derivative along the trajectories of~\eqref{eq:agents} with inputs~\eqref{eq:ContInput} given by $\dot{V}=x^T\bar{\mathcal{L}}x$.
\end{lemma}

From Theorem \ref{th:consensus} it can be seen that $\dot{V}$ is negative when the overall system is in disagreement and is equal to zero only when the corresponding states are in total agreement. In the latter case we also have $V=0$, see Lemma \ref{lm:LyapFn}. Different from consensus with single integrators, where the agents converge to a constant value, here it is only required that the difference between states of agents tends to zero, regardless of the particular response of the systems. As with many consensus algorithms, an estimate of the second smallest eigenvalue of the Laplacian matrix is required; this is the only global information needed by the agents.  Algorithms for distributed estimation of the second eigenvalue of the Laplacian have been presented in \cite{Aragues12}, \cite{Franceschelli09}. Readers are referred to these papers for details.


\section{Decentralized Periodic Event Triggered Consensus} \label{sec:four}
In this section we consider the case when agents use event-triggered communication strategies in order to reduce the frequency of transmissions. Every agent implements discrete-time models \eqref{eq:models} and samples its own state every $h$ time units in order to compute its local control input and check its local thresholds. It is assumed in this section that communication delays are negligible. We derive decentralized thresholds that depend only on local information and can be measured and applied in a decentralized way. Additionally, it is shown how the choice of the period $h$ affects the bounds on the disagreement between states of agents and how to select the threshold parameters to determine lower bounds on the inter-event time intervals. 

Define 
\begin{align}
	 e_i(t)=y_i(t)-x_i(t)  \label{eq:errei}
\end{align} 
and $x=\begin{bmatrix}
	x_1^T\ldots x_N^T
\end{bmatrix}^T$, $
y=\begin{bmatrix}
	y_1^T\ldots y_N^T
\end{bmatrix}^T$, $
e=\begin{bmatrix}
	e_1^T\ldots e_N^T
\end{bmatrix}^T$, 
where $y_i(t)$ represents the hypothetical continuous-time model
\begin{align}
	\dot{y}_i(t)=Ay_i(t) \ \ i=1...N. \label{eq:ContModel}
\end{align}
Definition \eqref{eq:ContModel} is only used for analysis of the event-triggered controller; it is not used for implementation of the event-triggered consensus algorithm.

Let us define the discretization errors
\begin{align}
	\left.
	\begin{array}{l l}
\breve{x}_i(t)=x_i(t_\mu)-x_i(t) \\
\breve{y}_i(t)=y_i(t_\mu)-y_i(t)  \\
\breve{e}_i(t)=e_i(t_\mu)-e_i(t)
\end{array} \right.  \label{eq:DiscErrors}
\end{align}
for $t\in[t_\mu,t_{\mu+1})$. Also define
$\breve{x}=\begin{bmatrix}
	\breve{x}_1^T\ldots \breve{x}_N^T
\end{bmatrix}^T$, $
\breve{y}=\begin{bmatrix}
	\breve{y}_1^T\ldots \breve{y}_N^T
\end{bmatrix}^T$, $
\breve{e}=\begin{bmatrix}
	\breve{e}_1^T\ldots \breve{e}_N^T
\end{bmatrix}^T$. 
Instead of $c$, the new coupling factor $c_1=2c$ is now used in the inputs \eqref{eq:inputs}.

\begin{theorem}\label{th:Decent}
 Assume the pair ($A,B$) is controllable and the communication graph is connected and undirected. Define $F$ in~\eqref{eq:F} and 
$c_1=2c$ where $c\geq 1/\lambda_2$. Then agents~\eqref{eq:agents} with inputs~\eqref{eq:inputs} and discrete-time models \eqref{eq:models} achieve a bounded consensus error where the difference between any two states is bounded by
\begin{align}
\left.
	\begin{array}{l l}
  \lim_{t\rightarrow\infty}\left\|x_i(t)-x_j(t)\right\|^2 \leq \frac{N\eta}{\beta\lambda_{\min}(P)} \label{eq:disagreement}
\end{array} \right.
\end{align}	
for $i,j=1,...,N$ and for any $t\geq 0$, if the local events are triggered when     
\begin{align}
  \delta_i>\sigma c_1 z^T_iPBB^TPz_i+\eta, \label{eq:DecCond}
\end{align}
where $0<\sigma<1$, $\beta=\frac{\lambda_{\min\neq 0}(-\bar{\mathcal{L}})}{\lambda_{\max}(\hat{\mathcal{L}})}>0$, 
\begin{align}
\left.
	\begin{array}{l l}
    \delta_i\!\!\!\!&=cN_i(b(N-1)+\frac{3N-1}{b})e^T_i(t_\mu)PBB^TPe_i(t_\mu)  \\
		&~~+c_1(1+2bN_i)\breve{z}^T_i(t^-_\mu)PBB^TP\breve{z}_i(t^-_\mu)  \\ 
		&~~+cN_i(\frac{N+1}{b}\!+\!3b(N\!-\!1)) \breve{e}^T_i(t^-_\mu)PBB^TP\breve{e}_i(t^-_\mu)  \\
	\end{array} \label{eq:kekz} \right.
\end{align}
$b>0$, 
\begin{align}
z_i=\sum_{j\in\mathcal{N}_i}\left(x_i(t_\mu)-y_j(t_\mu)\right).  \label{eq:zi}
\end{align}
and
\begin{align}
\breve{z}_i=\sum_{j\in\mathcal{N}_i}\left(\breve{x}_i-\breve{y}_j\right).  \label{eq:zicur}
\end{align}
 Furthermore, the agents do not exhibit Zeno behavior and the inter-event times $t_{k_i+1}-t_{k_i}$ for every agent $i=1,...,N$ are bounded as follows
\begin{align}
  0<h< t_{k_i+1}-t_{k_i}  \label{eq:tk}
\end{align}	
if 
\begin{align}
\left.
\begin{array}{l l}
  \eta\!>\big(2c_1N_i(b(N-1)+\frac{N}{b})\left\|E\right\|^2  + c_1(1+2bN_i)(2 + b_e\left\|E\right\|)^2 \big)\left\|PBB^TP\right\|\bar{z_i}^2 
\end{array} \label{eq:eta} \right.
\end{align}
where 
\begin{align}
  \bar{z}_i=\lambda_{\max} (\mathcal{L})\sqrt{\frac{V_M}{\lambda_{\min}(\hat{\mathcal{L})}}}    \label{eq:zibar} 
\end{align}
$V_M=\max \left\{V(0),\frac{N\eta}{\beta}\right\}$, $E=\int_0^h e^{A(h-s)}cBFds$, $b_e=\sqrt{N_iN(N-1)(\frac{b}{2}+\frac{1}{2b})}$, and $h=t_{\mu+1}-t_\mu$ is the discretization period.
\end{theorem}
\textit{Proof}. 
By implementing the coupling factor $c_1=2c$ in \eqref{eq:inputs}, we can write \eqref{eq:agents}-\eqref{eq:inputs} in compact form as follows:
\begin{align}
\left.
	\begin{array}{l l}
   \dot{x}\!=\bar{A}x+\bar{B}_\mathcal{D}x(t_\mu)+\bar{B}_\mathcal{A}y(t_\mu) 
			=(A_c+\bar{B})x+\bar{B}_\mathcal{A}e(t_\mu)+\bar{B}_1\breve{x} 
\end{array} \label{eq:dotx} \right.
\end{align}	
where $\bar{B}_\mathcal{D}=c_1\mathcal{D}\otimes BF$, $\bar{B}_\mathcal{A}=-c_1\mathcal{A}\otimes BF$, $\bar{B}_1=c_1\mathcal{L}\otimes BF$. The overall system dynamics \eqref{eq:dotx} is written in terms of the closed-loop state $x$ plus two error variables that are introduced in the system dynamics because of the sampled inputs and the event-based communication strategies. The error $e(t_\mu)$ is due to the fact that each agent only communicates its local state at some local event time instants and the error $\breve{x}$ is due to the fact that only discrete models of neighbors and periodic samples of the local state are used to compute the control input instead of continuous variables. In other words, the first error results from limiting the communication between agents and the second error results from reducing actuation time instants at each local node.

Consider the candidate Lyapunov function $V=x^{T}\hat{\mathcal{L}}x$ and evaluate the derivative along the trajectories of systems~\eqref{eq:agents} with inputs~\eqref{eq:inputs}. We can express $\dot{V}$ as follows: 
\begin{align}
\left.
	\begin{array}{l l}
	\dot{V}\!\!\!\!\!&=x^{T}\hat{\mathcal{L}}\left((A_c+\bar{B})x+\bar{B}_\mathcal{A}e(t_\mu)+\bar{B}_1\breve{x}\right)   
	+\left((A_c+\bar{B})x+\bar{B}_\mathcal{A}e(t_\mu)+\bar{B}_1\breve{x}\right)^{T}\hat{\mathcal{L}}x\\ 
	&=x^T\bar{\mathcal{L}}x + 2x^T\hat{\mathcal{L}}\bar{B}x + 2x^T\hat{\mathcal{L}}\bar{B}_\mathcal{A}e(t_\mu) 
	 + 2x^T\hat{\mathcal{L}}\bar{B}_1\breve{x}  \\
	&=x^T\bar{\mathcal{L}}x 
	+2\sum^{N}_{i=1}\Big[\sum_{k\in\mathcal{N}_i}(x_i-x_k)^TPBB^TP \\
	&~~~~ \Big(\!-c\sum_{j\in\mathcal{N}_i}(x_i-x_j) + c_1\sum_{j\in\mathcal{N}_i}e_j(t_\mu)  \\
	&~~~~-c_1\sum_{j\in\mathcal{N}_i}(\breve{x}_i-\breve{x}_j)\Big)\Big]   \\
	\end{array} \label{eq:Vdot} \right.
\end{align}
Eq. \eqref{eq:Vdot} can be further written in terms of the sampled state information $x_i(t_\mu)$, $x_j(t_\mu)$ and of the discretization errors $\breve{x}_i$, $\breve{x}_j$ defined in \eqref{eq:DiscErrors} as follows:	
\begin{align}
\left.
	\begin{array}{l l}
  \dot{V}\!\!\!\!\!\!&=x^T\bar{\mathcal{L}}x 
	+2\!\!\sum^{N}_{i=1}\Big[\!\sum_{k\in\mathcal{N}_i}(x_i(t_\mu)\!-\!x_k(t_\mu)\!-\!(\breve{x}_i-\breve{x}_k))^T \\
	&~~~~PBB^TP\Big(\!-\!c\sum_{j\in\mathcal{N}_i}(x_i(t_\mu)\!-\!x_j(t_\mu)-(\breve{x}_i\!-\!\breve{x}_j)) \\
	&~~~~+c_1\sum_{j\in\mathcal{N}_i}e_j(t_\mu)-c_1\sum_{j\in\mathcal{N}_i}(\breve{x}_i-\breve{x}_j)\Big)\Big] \\
	&=x^T\bar{\mathcal{L}}x +2\!\!\sum^{N}_{i=1}\Big[\!\sum_{k\in\mathcal{N}_i}(x_i(t_\mu)\!-\!x_k(t_\mu)\!)^T\\
	&~~~~PBB^TP\Big(\!-\!c\!\sum_{j\in\mathcal{N}_i}(x_i(t_\mu)\!-\!x_j(t_\mu))\! \\
	&~~~~  + c_1\sum_{j\in\mathcal{N}_i}e_j(t_\mu)\Big) + \xi_i \Big]
  \end{array} \label{eq:VdotII} \right.
\end{align}
where 
\begin{align}
\left.
	\begin{array}{l l}
\xi_i\!\!\!\!&=\sum_{k\in\mathcal{N}_i}(\breve{x}_i\!-\!\breve{x}_k\!)^TPBB^TP\Big(c\!\sum_{j\in\mathcal{N}_i}(\breve{x}_i\!-\!\breve{x}_j)\! 
	 - \!c_1\sum_{j\in\mathcal{N}_i}e_j(t_\mu)\Big)
		\end{array} \label{eq:Xi} \right.
\end{align}
By using the state error definition \eqref{eq:errei} and the definition of $z_i$ in \eqref{eq:zi} we can write \eqref{eq:VdotII} in the following form:
\begin{align}
\left.
	\begin{array}{l l}
	\dot{V}\!\!\!\!&=x^T\bar{\mathcal{L}}x +2\sum^{N}_{i=1}\Big[(z_i+\sum_{k\in\mathcal{N}_i}e_k(t_\mu))^T  \\
	&~~PBB^TP\Big(\!-\!c(z_i+\sum_{j\in\mathcal{N}_i}e_j(t_\mu))\! 
	+ \!c_1\sum_{j\in\mathcal{N}_i}e_j(t_\mu)\Big) + \xi_i \Big] \\
	&=x^T\bar{\mathcal{L}}x +2\sum^{N}_{i=1}\Big[-cz_i^TPBB^TPz_i   \\
	&~~+c\sum_{k\in\mathcal{N}_i}e_k^T(t_\mu)PBB^TP\sum_{j\in\mathcal{N}_i}e_j(t_\mu) 
	+ \xi_i \Big]
	\end{array} \label{eq:Vdot2} \right.
\end{align}
Using the inequality $\left\|x^{T}y\right\|\leq\frac{b}{2}x^{T}x+\frac{1}{2b}y^{T}y$, for $b>0$, we have that the following expression involving the second term in~\eqref{eq:Vdot2} holds
\begin{align}
\left.
	\begin{array}{l l}
	\sum_{k\in\mathcal{N}_i}e^T_k(t_\mu)PBB^TP\sum_{j\in\mathcal{N}_i}e_j(t_\mu)  \\
	\leq \left\|\sum_{k\in\mathcal{N}_i}e^T_kPBB^TP\sum_{j\in\mathcal{N}_i}e_j\right\|\\
	\leq \sum_{k\in\mathcal{N}_i}\sum_{j\in\mathcal{N}_i}\left\|e^T_kPBB^TPe_j\right\|\\
	\leq N_i(\frac{b}{2}+\frac{1}{2b})\sum_{j\in\mathcal{N}_i}e^T_jPBB^TPe_j
	\end{array} \label{eq:eiejbound} \right.
\end{align}
Since the communication graph is undirected we have the following properties 
\begin{align}
  \sum^{N}_{i=1}\sum_{j\in\mathcal{N}_i}e^T_jPBB^TPe_j=\sum^{N}_{i=1}\sum_{j\in\mathcal{N}_i}e^T_iPBB^TPe_i \label{eq:UndProp1}
\end{align}
and
\begin{align}
\left.
	\begin{array}{l l}
  \sum^{N}_{i=1}N_i\sum_{j\in\mathcal{N}_i}e^T_jPBB^TPe_j\\
	=\sum^{N}_{i=1}N_j\sum_{j\in\mathcal{N}_i}e^T_iPBB^TPe_i\\
	\leq \sum^{N}_{i=1}(N-1)\sum_{j\in\mathcal{N}_i}e^T_iPBB^TPe_i \label{eq:UndProp2}
\end{array} \right.
\end{align}
Then, using~\eqref{eq:UndProp2}, we can write the following expression which bounds the second term in \eqref{eq:Vdot2}
\begin{align}
\left.
	\begin{array}{l l}
	\sum^{N}_{i=1}N_i(\frac{b}{2}+\frac{1}{2b})\sum_{j\in\mathcal{N}_i}e^T_jPBB^TPe_j\\
	=\sum^{N}_{i=1}N_j(\frac{b}{2}+\frac{1}{2b})\sum_{j\in\mathcal{N}_i}e^T_iPBB^TPe_i\\
	\leq \sum^{N}_{i=1}N_i(N-1)\left(\frac{b}{2}+\frac{1}{2b}\right)e^T_iPBB^TPe_i.
	\end{array} \label{eq:Sumjk} \right.
\end{align}
Let us now analyze the term \eqref{eq:Xi}. Using the definitions \eqref{eq:DiscErrors} and \eqref{eq:zicur} the summation of terms $\xi_i$ can be written as follows 
\begin{align}
\left.
	\begin{array}{l l}
\sum^{N}_{i=1}\xi_i& \!\!\!\!\!=\sum^{N}_{i=1}\Big[c\breve{z}_i^TPBB^TP\breve{z}_i \\
&+ c_1\breve{z}_i^TPBB^TP\sum_{j\in\mathcal{N}_i}\breve{e}_j  \\
&+ c\sum_{k\in\mathcal{N}_i}\breve{e}_k^TPBB^TP\sum_{j\in\mathcal{N}_i}\breve{e}_j \\
&- c_1\breve{z}_i^TPBB^TP\sum_{j\in\mathcal{N}_i}e_j(t_\mu)  \\
&+ c_1\sum_{k\in\mathcal{N}_i}\breve{e}_k^TPBB^TP\sum_{j\in\mathcal{N}_i}e_j(t_\mu) \Big]
	\end{array} \label{eq:DiscrErrors} \right.
\end{align}
We again use the inequality $\left\|x^{T}y\right\|\leq\frac{b}{2}x^{T}x+\frac{1}{2b}y^{T}y$, for $b>0$, to obtain:
\begin{align}
\left.
	\begin{array}{l l}
\sum^{N}_{i=1}\xi_i & \!\!\!\!\!\leq \sum^{N}_{i=1}\Big[c(1+2bN_i)\breve{z}_i^TPBB^TP\breve{z}_i \\
&+ cN_i(\frac{N+1}{2b}+\frac{3b(N-1)}{2}) \breve{e}_i^TPBB^TP\breve{e}_i  \\
&+ \frac{c}{b}N_iN e_i^T(t_\mu)PBB^TPe_i(t_\mu)  \Big] \\
	\end{array} \label{eq:DiscrErrorsII} \right.
\end{align}
Using \eqref{eq:Sumjk} and \eqref{eq:DiscrErrorsII} we can bound \eqref{eq:Vdot2} as follows:
\begin{align}
\left.
	\begin{array}{l l}
	\dot{V}\!\!\!\!\!&\leq x^T\bar{\mathcal{L}}x +\sum^{N}_{i=1}\Big[-c_1z_i^TPBB^TPz_i   \\
	&~~+cN_i(b(N\!-\!1)+\frac{3N\!-\!1}{b})e^T_i(t_\mu)PBB^TPe_i(t_\mu)  \\
		&~~+c_1(1+2bN_i)\breve{z}^T_iPBB^TP\breve{z}_i  \\ 
		&~~+cN_i(\frac{N+1}{b}\!+\!3b(N\!-\!1)) \breve{e}^T_iPBB^TP\breve{e}_i  \Big]
	\end{array} \label{eq:Vdot3} \right. 
\end{align}
Note that we can compute errors and make communication decisions only at the sampling time instants $t_\mu$, that is, only when each agent measures its local state. Also note that the variables $\breve{z}_i$ and $\breve{e}_i$ are equal to zero at every time $t_\mu$ regardless if an event is triggered or not, i.e. $\breve{z}_i(t_\mu)=\breve{e}_i(t_\mu)=0$, for $\mu=0,1,2,...$. However, it is necessary to consider the non-zero value (before they are reset to zero) of these errors in order to evaluate \eqref{eq:Vdot3} and to decide if an event needs to be triggered or not. In other words, at every time instant $t_\mu$ we need to evaluate event thresholds using the discretization errors just before they are reset to zero, that is, just before the actuation update occurs at time $t_\mu$. 

\begin{figure}
	\begin{center}
	\begin{psfrags}
	\psfrag{ybjt}{$\breve{y}_j(t_\mu^-)$}
	\psfrag{xbit}{$\breve{x}_i(t_\mu^-)$}
		\includegraphics[width=8.4cm,height=6.5cm,trim=1cm .1cm .8cm .1cm]{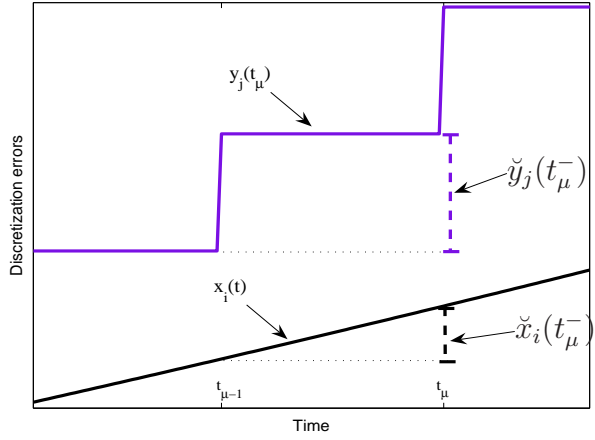}
	\end{psfrags}
	\caption{Computation of discretization errors based on current and previous measurements}
	\label{fig:DiscErr}
	\end{center}
\end{figure}

Let us denote the time just before the actuation update as $t_\mu^-$, so we check the discretization errors in \eqref{eq:Vdot3} using $\breve{z}_i(t_\mu^-)$ and $\breve{e}_i(t_\mu^-)$, i.e. the errors just before they are reset to zero. Thus, we can write 
\begin{align}
	\dot{V}&\leq x^T\bar{\mathcal{L}}x +\sum^{N}_{i=1}\big[-c_1z_i^TPBB^TPz_i +\delta_i \big]  \label{eq:Vdot4}
\end{align}
where $\breve{z}_i(t_\mu^-)=\sum_{j\in\mathcal{N}_i}(\breve{x}_i(t_\mu^-)-\breve{y}_j(t_\mu^-))$, $\breve{e}_i(t_\mu^-)=\breve{y}_i(t_\mu^-)-\breve{x}_i(t_\mu^-)$, and the variables $\breve{x}_i(t_\mu^-)$, $\breve{y}_i(t_\mu^-)$, and $\breve{y}_j(t_\mu^-)$, for $j\in \mathcal{N}_i$ can be obtained locally (see Fig. \ref{fig:DiscErr}) at every time $t_\mu$ using the current and previous samples of the local state and the local models as follows 
\begin{align}
\left.
	\begin{array}{l l}
   \breve{x}_i(t_\mu^-)=x_i(t_{\mu-1})-x_i(t_{\mu})  \\
   \breve{y}_j(t_\mu^-)=y_j(t_{\mu-1})-y_j(t_{\mu}).  
	\end{array} \label{eq:DExiyj} \right. 
\end{align}
The variables $\breve{x}_i$ and $\breve{y}_j$, $j=i,j\in \mathcal{N}_i$, cannot be measured continuously but only at times $t_\mu$ as it is illustrated in Fig. \ref{fig:DiscErr}. 

When threshold~\eqref{eq:DecCond} holds, say, at time $t_{k_i}$ then the error resets to zero, that is, $e_i(t_{k_i})=0$, since $y_i(t_{k_i})=x_i(t_{k_i})$. Also, since $t_{k_i}=t_\mu$ for some non-negative integer $\mu$ (recall that communication events can only occur at the sampling times $t_\mu$), then, the control input is updated and the discretization errors are equal to zero, that is, $\breve{z}_i(t_\mu)=\breve{e}_i(t_\mu)=0$. This means that $\delta_i(t_{k_i})=0$, then the following holds  
\begin{align}
 \delta_i \leq \sigma c_1z_i^TPBB^TPz_i + \eta. 
\end{align}
Consequently,
\begin{align}
  \dot{V}&\leq x^{T}\bar{\mathcal{L}}x+\sum^{N}_{i=1}\left[(\sigma-1)c_1z_i^{T}PBB^TPz_i +\eta \right]
\end{align}
Since $\sigma-1<0$ we have that $(\sigma-1)c_1z_i^{T}PBB^TPz_i\leq 0$ for $i=1,...,N$ and	
\begin{align}
  \dot{V}&\leq x^{T}\bar{\mathcal{L}}x + N\eta
\end{align}
We use the fact that $\hat{\mathcal{L}}$ is positive semi-definite and it has $n$ zero eigenvalues with corresponding eigenvectors $v_1\ldots v_n$. Let $x=x_1+x_2$ such that $\left\langle x^T_1,x_2\right\rangle =0$ and $\hat{\mathcal{L}}x_1=0$, that is, $x_1$ belongs to the subspace spanned by $v_1\ldots v_n$. Consider
\begin{align}
  x^T\hat{\mathcal{L}}x=(x_1+x_2)^T\hat{\mathcal{L}}(x_1+x_2)=x_2^T\hat{\mathcal{L}}x_2
\end{align}
also
\begin{align}
  x_2^T\hat{\mathcal{L}}x_2\leq \lambda_{\max}(\hat{\mathcal{L}})x_2^Tx_2.
\end{align}
From Theorem~\ref{th:consensus} we know that $-\bar{\mathcal{L}}$ is positive semi-definite with $n$ zero eigenvalues and $\bar{\mathcal{L}}x_1=0$. We can see that
\begin{align}
x^T\bar{\mathcal{L}}x=(x_1+x_2)^T\bar{\mathcal{L}}(x_1+x_2)=x_2^T\bar{\mathcal{L}}x_2
\end{align}
and
\begin{align}
  x_2^T(-\bar{\mathcal{L}})x_2\geq \lambda_{\min\neq 0}(-\bar{\mathcal{L}})x_2^Tx_2
\end{align}
where $\lambda_{\min\neq 0}(-\bar{\mathcal{L}})>0$ is the minimum eigenvalue of $-\bar{\mathcal{L}}$ other than zero. Combining the expressions above we have
\begin{align}
  x_2^T\bar{\mathcal{L}}x_2\leq -\frac{\lambda_{\min\neq 0}(-\bar{\mathcal{L}})}{\lambda_{\max}(\hat{\mathcal{L}})}x_2^T\hat{\mathcal{L}}x_2.
\end{align}
Then $\dot{V}$ can be bounded as follows:
\begin{align}
  \dot{V}\leq -\beta x^{T}\hat{\mathcal{L}}x+N\eta=-\beta V+N\eta. \label{eq:Vdiffineq}
\end{align}
Solving~\eqref{eq:Vdiffineq} we have that
\begin{align}
\left.
	\begin{array}{l l}
  V(t)&\leq e^{-\beta t}V(0)+ N\eta\int^t_{0}e^{-\beta (t-\tau)}d\tau\\
   &\leq \left(V(0)-\frac{N\eta}{\beta}\right)e^{-\beta t}+\frac{N\eta}{\beta}.
\end{array} \label{eq:Voft} \right.
\end{align}
Expression~\eqref{eq:Voft} represents a bound on the consensus states as a function of the initial separation of the agents $V(0)=x(0)^T\hat{\mathcal{L}}x(0)$. We can express $V(t)=\frac{1}{2}\sum^N_{i=1}\sum_{j\in\mathcal{N}_i}(x_i-x_j)^TP(x_i-x_j)$ and a direct bound on the difference between any two states $i,j$ can be obtained as follows. Since the graph is undirected the term $(x_i-x_j)^TP(x_i-x_j)$ appears twice in the summation $V(t)$, then we can write
\begin{align}
\left.
	\begin{array}{l l}
  \lambda_{\min}(P)\left\|x_i-x_j\right\|^2\!\!\!\! &\leq (x_i-x_j)^TP(x_i-x_j)\\
	&\leq \!\!\left(V(0)-\frac{N\eta}{\beta}\right)e^{-\beta t}+\frac{N\eta}{\beta}.
\end{array} \right.
\end{align}	
Finally, the difference between any two states can be bounded as in \eqref{eq:disagreement}.

We will now prove that the inter-event times are lower bounded by the sampling period as in \eqref{eq:tk}. 

First note that 
\begin{align}
\left.
	\begin{array}{l l}
   \left\|z_i\right\|\!\!\!&=\left\| \sum_{j\in\mathcal{N}_i}(x_i(t_\mu) - x_j(t_\mu) - e_j(t_\mu))\right\|   \\
 &\leq \left\|\sum_{j\in\mathcal{N}_i}(x_i(t_\mu) - x_j(t_\mu))\right\| + \left\|\sum_{j\in\mathcal{N}_i}e_j(t_\mu)\right\|    \\
&\leq \left\| \mathcal{L}_n x\right\| + \left\| \mathcal{A}_n e\right\|
\end{array} \right.  \label{eq:zib}
\end{align}	
where $\mathcal{L}_n=\mathcal{L}\otimes I_n$ and $\mathcal{A}_n=\mathcal{A}\otimes I_n$. For the first term in \eqref{eq:zib} we have that the following holds
\begin{align}
\left.
	\begin{array}{l l}
  \left\| \mathcal{L}_n x\right\|\!\!\!&=\sqrt{x^T\mathcal{L}_n \mathcal{L}_n x} \leq \lambda_{\max}(\mathcal{L}_n)\sqrt{x^Tx}   \\
	  &\leq  \lambda_{\max}(\mathcal{L}_n)\sqrt{\frac{V(t_\mu)}{\lambda_{\min}(\hat{\mathcal{L}})}} \\
		&\leq  \lambda_{\max}(\mathcal{L})\sqrt{\frac{V_M}{\lambda_{\min}(\hat{\mathcal{L}})}}
\end{array} \right.  \label{eq:Lapnx}
\end{align}	
Te second term in \eqref{eq:zib} can be written as follows
\begin{align}
\left.
	\begin{array}{l l}
  \left\| \mathcal{A}_n e\right\|\!\!\!&=\sqrt{e^T\mathcal{A}_n \mathcal{A}_n e}  \\
	&= \sqrt{\sum_{i=1}^N \sum_{k\in\mathcal{N}_i}e_k^T \sum_{j\in\mathcal{N}_i}e_j}  \\
	&\leq \sqrt{\sum_{i=1}^N N_i(N-1)(\frac{b}{2}+\frac{1}{2b})}e_i^Te_i  \\
	&\leq b_e \left\|e_i\right\|
\end{array} \right.  \label{eq:Ane}
\end{align}	
Then we can write \eqref{eq:zib} as follows
\begin{align}
\left.
	\begin{array}{l l}
   \left\|z_i\right\|\!\!\!&\leq \lambda_{\max}(\mathcal{L})\sqrt{\frac{V_M}{\lambda_{\min}(\hat{\mathcal{L}})}} + b_e \left\|e_i\right\|
\end{array} \right.  \label{eq:zib2}
\end{align}	
for any time $t_\mu$. In particular, at the local event times $t_{k_i}$ we have that $e_i(t_{k_i})=0$ and $\left\|z_i(t_{k_i})\right\|\leq \bar{z}_i$.

In order to prove \eqref{eq:tk} we need to guarantee that \eqref{eq:kekz} does not grow from zero at time $t_{k_i}$ to $\eta$ at time $t_{k_i}+h$ and will not trigger an event at $t_{k_i}+h$. In other words, for a given sampling period $h$, we need to find the value of $\eta$ that guarantees \eqref{eq:tk}. 

Consider the error dynamics of agent $i$ from $t_{k_i}$ to $t_{k_i}+h$ as follows:
\begin{align}
 \left.
	\begin{array}{l l}
 e_i(t_{k_i}+h)\!\!\!&=y_i(t_{k_i}+h) - x_i(t_{k_i}+h) \\ 
	&=Gy_i(t_{k_i})-Gx_i(t_{k_i})-Ez_i(t_{k_i}) \\
	&=Ge_i(t_{k_i})-Ez_i(t_{k_i}). 
	\end{array} \right.  \label{eq:erroriDyn}
\end{align}	
Note that $e_i(t_{k_i})=0$ since the error is reset at the event time $t_{k_i}$, then we have that
\begin{align}
 e_i(t_{k_i}+h)&=-Ez_i(t_{k_i}).   \label{eq:erroriDyn2}
\end{align}	
The effect of the discretization errors is considered before they are reset to zero at time $t_{k_i}+h$ because of the actuation update, so we have
\begin{align}
 \left.
	\begin{array}{l l}
 \breve{e}_i(t_{k_i}+h)\!\!\!&=e_i(t_{k_i}) - e_i(t_{k_i}+h) \\ 
	&=- e_i(t_{k_i}+h)=Ez_i(t_{k_i}). 
	\end{array} \right.  \label{eq:CurerroriDyn}
\end{align}	
Now, let us analyze the term $\breve{z}_i$. Consider
\begin{align}
 \left.
	\begin{array}{l l}
 \left\|\breve{z}_i(t_{k_i}+h)\right\|\!   
 &=\left\| \sum_{j\in\mathcal{N}_i}(\breve{x}_i(t_{k_i}+h) - \breve{x}_j(t_{k_i}+h) - \breve{e}_j(t_{k_i}+h))\right\| \\
&\leq \left\| \mathcal{L}_n\breve{x}(t_{k_i}+h)\right\| + b_e\left\|\breve{e}_i(t_{k_i}+h)\right\|
	\end{array} \right.  \label{eq:Curzi}
\end{align}	
Note that
\begin{align}
 \left.
	\begin{array}{l l}
 \left\|\mathcal{L}_n \breve{x} (t_{k_i}+h)\right\|  
  &=\left\|\mathcal{L}_n x(t_{k_i}) - \mathcal{L}_n x(t_{k_i}+h)\right\|  \\
  &\leq \frac{ \lambda_{\max}(\mathcal{L})}{\sqrt{\lambda_{\min}(\hat{\mathcal{L}})}} \big(\sqrt{V(t_{k_i})}+\sqrt{V(t_{k_i}+h)}\big)  \\
	&\leq \frac{ 2\lambda_{\max}(\mathcal{L})}{\sqrt{\lambda_{\min}(\hat{\mathcal{L}})}} \sqrt{V_M}
	\end{array} \right.  \label{eq:Lapy}
\end{align}	
We use \eqref{eq:erroriDyn2}, \eqref{eq:CurerroriDyn}, and \eqref{eq:Lapy} to analyze the growth of the term $\delta_i$ from $t_{k_i}$ to $t_{k_i}+h$
\begin{align}
\left.
	\begin{array}{l l}
  \delta_i(t_{k_i}+h) 
	&\leq  2c_1N_i(b(N\!-\!1)+\frac{N}{b})\left\|E\right\|^2\left\|PBB^TP\right\|z_i^2(t_{k_i}) \\
	 &  ~~+c_1(1+2bN_i)\left\|PBB^TP\right\|  
		\Big(\frac{ 2\lambda_{\max}(\mathcal{L})}{\sqrt{\lambda_{\min}(\hat{\mathcal{L}})}}\sqrt{V_M} +b_e\left\|E\right\|z_i(t_{k_i})\Big)^2 
	\end{array} \right.  \label{eq:errorgrow}
\end{align}
Then, by the selection of $\eta$ in \eqref{eq:eta} we can guarantee that after an event instant at $t_{k_i}$ the term $\delta_i$ cannot grow from zero to $\eta$ and threshold \eqref{eq:DecCond} is not triggered. Then, the inter-event times are greater than $h>0$ and, obviously, we can guarantee that Zeno behavior does not occur at any node. $\bullet$

\textit{Remark}. The design parameter $\eta$ provides a tradeoff between performance as measured by the consensus error \eqref{eq:disagreement} and the bound on the inter-event intervals \eqref{eq:tk}. Also note that the variables used to compute the threshold~\eqref{eq:DecCond}, which define the events at node $i$, are available locally. Concerning global information we only need an estimate of the second eigenvalue of the Laplacian, as it was mentioned earlier.

\textit{Remark}. Note that the triggering of a local event by agent $i$, i.e. the transmission of a measurement $x_i(t_{k_i})$, does not change the local control input $u_i$ beyond the local actuation update, since $u_i$ is not a function of $y_i$. The transmission of a new measurement updates the control inputs of neighbor agents and for the local agent, it only resets its local state error $e_i$ and $\delta_i$. 

\textit{Remark}. The results provided in this paper hold for agents described by general linear dynamics. The common cases of single integrators and double integrators are particular cases covered by this framework. The single integrator is modeled as a ZOH and the double integrator is modeled similar to \cite{Seyboth13}, that is, velocity as a ZOH and position as a first-order-hold model. 
 
\textit{Remark}. The event-triggered consensus algorithm  provides a bound on the difference between any two states not only at the sampling time instants but for any $t\geq 0$. Also note that $\eta$ can be made arbitrarily small (thus decreasing the bound \eqref{eq:disagreement}) and still satisfy \eqref{eq:tk} by choosing a small discretization period $h$ since we have that both $E\rightarrow 0$ and $I-G\rightarrow 0$ (therefore $\delta_i(t_{k_i}+h) \rightarrow 0$) as $h\rightarrow 0$.

\section{Decentralized Periodic Event Triggered Consensus With Communication Delays} \label{sec:five}

In this section we consider the presence of time-varying but bounded communication delays. Since the measurement updates will be delayed, the neighbors of an agent $i$ will have a version of agent $i$'s model state that is different than agent $i$'s version. It is necessary to distinguish between the model state as seen by the local agent, itself, and as seen by its neighbors. Define the dynamics and update law of the model state of agent $i$ as seen by agent $i$ as
\begin{align}
  y_{ii}(t_{\mu+1})=Gy_{ii}(t_\mu),  \;  y_{ii}(t_{k_i})=x_{i}(t_{k_i}). \label{eq:yii}
\end{align}
The measurement $x_{i}(t_{k_i})$ is transmitted by agent $i$ at time $t_{k_i}$ and will arrive at agents $j$, $j\in \mathcal{N}_i$, at time $t_{k_i}+d_i(t_{k_i})$. For a given update instant all receiving agents experience the same delay $d_i(t_{k_i})$. However, this is not a constraint and the communication delays can be generalized so, for the same update instant $t_{k_i}$, the neighbors are updated at different time instants. This is described at the end of this section. 

Assume without loss of generality that $d_i(t_{k_i})=hp_i(t_{k_i})$ and $p_i\geq 1$ is an integer, that is, the delay is an integer multiple of the sampling period $h$. If the delay is not an integer multiple of the sampling period the receiving agent uses the corresponding delayed measurement at the next sampling period to update the corresponding model which effectively makes the delay to be an integer multiple of $h$. 

Define the dynamics of the model state of agent $i$ as seen by agent $j$, $j\in \mathcal{N}_i$, as
\begin{align}
    y_{ij}(t_{\mu+1})=Gy_{ij}(t_\mu),  \label{eq:yij}
\end{align}
The states of these models are updated when a delayed measurement of agent $i$ is received by agent $j$, $j\in \mathcal{N}_i$. The state measurement $x_i(t_{k_i})$ transmitted by agent $i$ at time instant $t_{k_i}$ is received by agent $j$, $j\in \mathcal{N}_i$, at time instant $t_{k_i}+hp_i(t_{k_i})$. Let us define the update law of the model state of agent $i$ as seen by agent $j$, $j\in \mathcal{N}_i$, as
\begin{align}
		y_{ij}(t_{k_i}+hp_i(t_{k_i}))=f_d(x_{i}(t_{k_i}),p_i(t_{k_i})) \label{eq:yijupd}
\end{align}
where $t_{k_i}$ represents the update instants triggered by agent $i$ and $hp_i(t_{k_i})$ represents the communication delay associated to the triggering instant $t_{k_i}$.

Define a positive and constant upper bound on the communication delays by $d=ph<t_{k_i+1}-t_{k_i}$, that is, $d_i(t_{k_i})\leq d$ for any triggering instant $t_{k_i}$ and for $i=1,...,N$. Later in this section we will define the design parameters that bound the inter-event times as a function of the delay $d$. Note that $p \geq p_i(t_{k_i})$, for $i=1,...,N$. Assume that the current delay, $p_i(t_{k_i})$, is known to the receiving agents, for instance, by applying time-stamping techniques. 

Since both, $y_{ii}$ and $y_{ij}$, use the same state matrix to compute their response between their corresponding update instants, then we define
\begin{align}
   f_d(x_{i}(t_{k_i}),hp_i(t_{k_i}))\triangleq G^{p_i(t_{k_i})}x_{i}(t_{k_i}) \label{eq:yij_prop}
\end{align}
that is, the delayed measurement is propagated forward in time and the result is used to update the state of the model as shown in \eqref{eq:yijupd}. By definition, we have that the following \textit{local} triggering event will occur at time $t_{k_i+1}>t_{k_i}+d$, this means that $y_{ii}(t_{k_i}+d_i(t_{k_i}))=G^{p_i(t_{k_i})}x_{i}(t_{k_i})$ because no other local event has been triggered since time instant $t_{k_i}$. Therefore, we have that $y_{ii}(t_\mu)\neq y_{ij}(t_\mu)$ for $t_\mu \in [t_{k_i},t_{k_i}+hp_i(t_{k_i}))$ and $y_{ii}(t_\mu)=y_{ij}(t_\mu)$ for $t\in [t_{k_i}+hp_i(t_{k_i}),t_{k_i+1})$.  

\begin{figure}
	\begin{center}
		\includegraphics[width=8.4cm,height=6.5cm,trim=1cm .2cm .8cm .1cm]{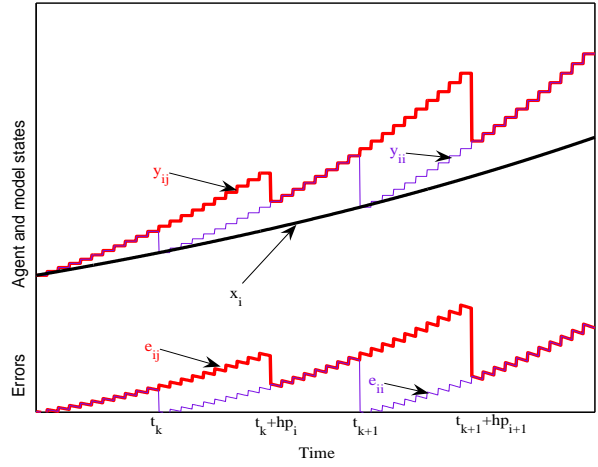}
	\caption{Relation between state $x_i$, model states $y_{ii}$, $y_{ij}$, and corresponding errors $e_{ii}$, $e_{ij}$.}
	\label{fig:FigureModErr}
	\end{center}
\end{figure}

Define the state errors
\begin{align}
  e_{ii}(t_\mu)=y_{ii}(t_\mu)-x_i(t_\mu), \label{eq:eii}  \\ 
  e_{ij}(t_\mu)=y_{ij}(t_\mu)-x_i(t_\mu).  \label{eq:eij}
\end{align}
Note that $e_{ii}(t_{k_i})=0$ and $e_{ij}(t_\mu)=e_{ii}(t_\mu)$, for $t_\mu\in [t_{k_i}+hp_i(t_{k_i}),t_{k_i+1})$. These relations are pictured in Fig. \ref{fig:FigureModErr}.

Also define $\zeta(t_\mu)=\begin{bmatrix}
	e_{11}^T(t_\mu)\ldots e_{NN}^T(t_\mu)
\end{bmatrix}^T\in \mathbb{R}^{nN}$ and $\zeta_d(t_\mu)=\begin{bmatrix}
	e_{1j_1}^T(t_\mu)\ldots e_{Nj_N}^T(t_\mu)
\end{bmatrix}^T\in \mathbb{R}^{nN}$, where the components $e_{ij_i}$ in $\zeta_d$ represent the errors defined in \eqref{eq:eij}, that is, the error of agent $i$ as seen by its neighbors $j_i$, $j_i\in \mathcal{N}_i$.

Define the discretization errors
$\breve{y}_{ii}(t)=y_{ii}(t_\mu)-y_{ii}(t)$, $\breve{y}_{ij}(t)=y_{ij}(t_\mu)-y_{ij}(t)$, $\breve{e}_{ii}(t)=e_{ii}(t_\mu)-e_{ii}(t)$, $\breve{e}_{ij}(t)=e_{ij}(t_\mu)-e_{ij}(t)$, 
$\breve{y}=\begin{bmatrix}
	\breve{y}_{11}^T\ldots \breve{y}_{NN}^T
\end{bmatrix}^T$, $\breve{y}_d=\begin{bmatrix}
	\breve{y}_{1j_1}^T\ldots \breve{y}_{Nj_N}^T
\end{bmatrix}^T$, $\breve{\zeta}=\begin{bmatrix}
	\breve{e}_{11}^T\ldots \breve{e}_{NN}^T
\end{bmatrix}^T$, $\breve{\zeta}_d=\begin{bmatrix}
	\breve{e}_{1j_1}^T\ldots \breve{e}_{Nj_N}^T
\end{bmatrix}^T$.

The dynamics of every agent in \eqref{eq:agents} with communication delays captured by the new input definitions (which are functions of delayed model states $y_{ji}(t_\mu)$): 
\begin{align}
	u_i(t_\mu)=c_1F\sum_{j\in\mathcal{N}_i}(x_i(t_\mu)-y_{ji}(t_\mu)),\ \ i=1...N, \label{eq:newinputs}
\end{align}
can be written in compact form as follows:
\begin{align}
   \dot{x}=(A_c+\bar{B})x+\bar{B}_\mathcal{A}\zeta_d(t_\mu)+\bar{B}_1\breve{x} 
\label{eq:dotxdel} 
\end{align}	
where the coupling strength $c_1=2c$ has been used.

Note that if we follow the same analysis as in Theorem \ref{th:Decent} we can arrive at an expression similar to \eqref{eq:Vdot3} involving the local control inputs $z_i$, the local discretization errors, and the delayed errors $e_{ij}(t_\mu)$. The latter will create a major difficulty in designing the local events since the local agent $i$ does not have access to the errors $e_{ij}(t_\mu)$ as seen by its neighbors. More importantly, the local agent is not able to reset the error $e_{ij}(t_\mu)$ but only the local error $e_{ii}(t_\mu)$. 

The following theorem provides a method to design local events in the presence of communication delays and using the local state errors $e_{ii}(t_\mu)$ to locally evaluate the events.

\begin{theorem}\label{th:DecentDelay}
 Assume the pair ($A,B$) is controllable and the communication graph is connected and undirected. Define $F$ in~\eqref{eq:F} and 
$c_1=2c$ where $c\geq 1/\lambda_2$. Then agents~\eqref{eq:agents} with inputs~\eqref{eq:newinputs} achieve, in the presence of communication delays $d_i<d$, a bounded consensus error where the difference between any two states is bounded by \eqref{eq:disagreement} for $i,j=1,...,N$, if the events are triggered when     
\begin{align}
  \delta_i>\sigma c_1z^T_iPBB^TPz_i+\eta  \label{eq:CondDelay}
\end{align}
where $0<\sigma<1$,  
\begin{align}
  \delta_i&= c_1(1+b)^2 \breve{z}^T_i(t_\mu^-)PBB^TP\breve{z}_i(t_\mu^-) + \delta_{di} + \breve{\delta}_i. \label{eq:deltadel}
\end{align}
\begin{align}
z_i=\sum_{j\in\mathcal{N}_i}\left(x_i(t_\mu)-y_{ji}(t_\mu)\right),  \label{eq:zidelay}
\end{align}
\begin{align}
\breve{z}_i=\sum_{j\in\mathcal{N}_i}\left(\breve{x}_i-\breve{y}_{ji}\right),  \label{eq:curzidelay}
\end{align}
The terms $\delta_{di}$ and $\breve{\delta}_i$ are given by
\begin{align}
\left.
\begin{array}{l l}
  \delta_{di}\!=c_1(1+\frac{1}{b})\bar{\lambda} \Big[e_{ii}^T(t_\mu)(G^p)^TG^pe_{ii}(t_\mu) 
	 + 2\left\|G^pe_{ii}(t_\mu)\right\|\Upsilon \bar{z}_i + (\Upsilon \bar{z}_i)^2  \Big], \label{eq:deltadd}
\end{array} \right.
\end{align} 
\begin{align}
\left.
\begin{array}{l l}
  \breve{\delta}_{i}\!\!\!\!&=c_1(1+b)(1+\frac{1}{b})\bar{\lambda} \Big[e_{ii}^T(t_\mu)((G-I)G^{p-1})^T
	 (G-I)G^{p-1}e_{ii}(t_\mu) \\
	&~~ + 2\left\|(G\!-\!I)G^{p-1}e_{ii}(t_\mu)\right\|  
  (\left\|G\!-\!I\right\|\Upsilon_h \!+\! \left\|E\right\| )\bar{z}_i \\
	&~~ + \big((\left\|G-I\right\|\Upsilon_h +\left\|E\right\| )\bar{z}_i\big)^2  \Big], \label{eq:deltacur}
\end{array} \right.
\end{align}
where $b>0$, $\bar{\lambda}=\lambda_{max}(\mathcal{A}^2\otimes PBB^TP)$, $\bar{z}_i=\frac{\lambda_{\max}(\mathcal{L})}{(1-b_e\Upsilon)}\sqrt{\frac{V_M}{\lambda_{\min}(\hat{\mathcal{L}})}}$, $b_e=\sqrt{N_iN(N-1)(\frac{b}{2}+\frac{1}{2b})}$, $\Upsilon=\int_0^d\left\|e^{A(d-s)}cBF\right\|ds<\frac{1}{b_e}$, $\Upsilon_h=\int_0^{d-h}\left\|e^{A(d-h-s)}cBF\right\|ds$, $E=\int_0^he^{A(h-s)}cBFds$.

Furthermore, the agents do not exhibit Zeno behavior and the inter-event times $t_{k_i+1}-t_{k_i}$ for every agent $i=1,...,N$ are bounded by $d>0$, that is
\begin{align}
  0<d<t_{k_i+1}-t_{k_i}.  \label{eq:tkdelay}
\end{align}	
if 
\begin{align}
\left.
\begin{array}{l l}
  \eta\!\!\!\!&>c_1(1+b)^2\left\|PBB^TP\right\| 
	\frac{\lambda_{max}^2(\mathcal{L})}{\lambda_{min}(\hat{\mathcal{L}})}V_M \Big(\frac{1}{1-b_e\Upsilon}+\frac{1}{1-b_e\Upsilon_h}  \Big)^2     \\
	&~~+c_1(1+\frac{1}{b})\big(\left\|G^p\right\|+1\big)^2\Upsilon^2 \bar{z}_i^2  \\
	&~~+c_1(1\!+\!b)(1\!+\!\frac{1}{b})\Big(\left\|(G\!-\!I)G^{p-1}\right\|\Upsilon  
	 +\left\|(G\!-\!I)\right\|\Upsilon_h\! + \!\left\|E\right\| \Big)^2 \bar{z}_i^2 	
\end{array} \right.  \label{eq:etadel}
\end{align}
\end{theorem}
\textit{Proof}. 
Consider the candidate Lyapunov function $V=x^{T}\hat{\mathcal{L}}x$ and evaluate the derivative along the trajectories of 
systems~\eqref{eq:agents} with inputs~\eqref{eq:newinputs}. 
\begin{align}
\left.
	\begin{array}{l l}
	\dot{V}\!\!\!\!&=x^{T}\hat{\mathcal{L}}\left((A_c+\bar{B})x+\bar{B}_\mathcal{A}\zeta_d(t_\mu)+\bar{B}_1\breve{x}\right)   \\
	&~~+\left((A_c+\bar{B})x+\bar{B}_\mathcal{A}\zeta_d(t_\mu)+\bar{B}_1\breve{x}\right)^{T}\hat{\mathcal{L}}x. 
	\end{array} \label{eq:Vdotdelay} \right.
\end{align}
Following similar steps to \eqref{eq:Vdot}-\eqref{eq:Vdot2} in the proof of Theorem \ref{th:Decent} and using the inequality $\left\|x^{T}y\right\|\leq\frac{b}{2}x^{T}x+\frac{1}{2b}y^{T}y$, for $b>0$,  we can write
\begin{align}
\left.
	\begin{array}{l l}
	\dot{V}\!\!\!\! &\leq x^T\bar{\mathcal{L}}x +2\sum^{N}_{i=1}\Big[-cz_i^TPBB^TPz_i   \\
	&~~+c(1\!+\!\frac{1}{b})\sum_{k\in\mathcal{N}_i}e_{ki}^T(t_\mu)  
	PBB^TP\sum_{j\in\mathcal{N}_i}e_{ji}(t_\mu) \\
	&~~+c(1\!+\!b)\sum_{k\in\mathcal{N}_i}(\breve{x}_i\!-\!\breve{x}_k\!)^T   
	PBB^TP\!\sum_{j\in\mathcal{N}_i}(\breve{x}_i\!-\!\breve{x}_j)\! \Big] 
	\end{array} \label{eq:Vdotdelay2} \right.
\end{align}
where $e_{ji}(t_\mu)$ represents the state error, at time instants $t_\mu$, of agent $j$ as seen by agent $i$, for $j\in\mathcal{N}_i$. We can write \eqref{eq:Vdotdelay2} as follows
\begin{align}
\left.
	\begin{array}{l l}
	\dot{V}\!\!\!\! &\leq x^T\bar{\mathcal{L}}x +\sum^{N}_{i=1}\Big[-c_1z_i^TPBB^TPz_i   \\
	&~~+c_1(1\!+\!\frac{1}{b})\sum_{k\in\mathcal{N}_i}e_{ki}^T(t_\mu)  
  PBB^TP\sum_{j\in\mathcal{N}_i}e_{ji}(t_\mu) \\
	&~~+c_1(1\!+\!b)\Big(\breve{z}_i^TPBB^TP\breve{z}_i  
   +2\breve{z}_i^TPBB^TP\sum_{j\in\mathcal{N}_i}\breve{e}_{ji}\\
	&~~ + \sum_{k\in\mathcal{N}_i}\breve{e}_{ki}^TPBB^TP\sum_{j\in\mathcal{N}_i}\breve{e}_{ji} \Big) \Big]. 
	\end{array} \label{eq:Vdotdelay2a} \right.
\end{align}
Using the inequality $\left\|x^{T}y\right\|\leq\frac{b}{2}x^{T}x+\frac{1}{2b}y^{T}y$, for $b>0$, once again, we obtain
\begin{align}
\left.
	\begin{array}{l l}
	\dot{V}\!\!\!\! &\leq x^T\bar{\mathcal{L}}x +\sum^{N}_{i=1}\Big[-c_1z_i^TPBB^TPz_i   \\
  &~~+c_1(1\!+\!b)^2\breve{z}_i^TPBB^TP\breve{z}_i  \\
	&~~+c_1(1\!+\!\frac{1}{b})\sum_{k\in\mathcal{N}_i}e_{ki}^T(t_\mu)PBB^TP\sum_{j\in\mathcal{N}_i}e_{ji}(t_\mu) \\
	&~~+c_1(1\!+\!b)(1\!+\!\frac{1}{b}) \sum_{k\in\mathcal{N}_i}\breve{e}_{ki}^TPBB^TP\sum_{j\in\mathcal{N}_i}\breve{e}_{ji}  \Big] \\ 
  &=x^T\bar{\mathcal{L}}x +\sum^{N}_{i=1}\Big[-c_1z_i^TPBB^TPz_i   \\
  &~~+c_1(1\!+\!b)^2\breve{z}_i^TPBB^TP\breve{z}_i \Big] \\
	&~~+c_1(1\!+\!\frac{1}{b})\zeta_d^T(t_\mu)(\mathcal{A}\otimes B^TP)^T(\mathcal{A}\otimes B^TP)\zeta_d^T(t_\mu)  \\
	&~~+c_1(1\!+\!b)(1\!+\!\frac{1}{b})\breve{\zeta}_d^T(\mathcal{A}\otimes B^TP)^T(\mathcal{A}\otimes B^TP)\breve{\zeta}_d^T.
	\end{array} \label{eq:Vdotdelay2b} \right.
\end{align}
The variables $z_i$ and $\breve{z}_i$ can be computed locally by every node. Let us then focus on the delayed error terms $\zeta_d^T(t_\mu)$ and $\breve{\zeta}_d^T$.
Define $\nu_i(t_\mu)=y_{ij}(t_\mu)-y_{ii}(t_\mu)$, then we have that
\begin{align}
  \left.
	\begin{array}{l l}
 e_{ij}(t_\mu)\!\!\!&=y_{ij}(t_\mu)-x_i(t_\mu) \\
	   &=y_{ij}(t_\mu)-(y_{ii}(t_\mu)-e_{ii}(t_\mu)) \\
		&=\nu_i(t_\mu)+e_{ii}(t_\mu).  
			\end{array}  \right. 	 \label{eq:errornu}
\end{align}	
For the term $\nu_i(t_\mu)$ the following holds:  
\begin{align}{\nu}_i(t_{\mu+1})=\left\{
	\begin{array}{l l}
     G\nu_i(t_\mu), \ \ t_\mu\in [t_{k_i},t_{k_i}+d_i(t_{k_i})) \\
	 0, \ \ \ \ \ \ \ \ \ \ t_\mu\in [t_{k_i}+d_i(t_{k_i}), t_{k_i+1})
 \end{array}  \right.  \label{eq:nudyn}
\end{align}	
with $\nu_i(t_{k_i})=e_{ii}(t_{k_i}^-)$. The notation $t_{k_i}^-$ represents the event time $t_{k_i}$ but just before the local error is reset to zero. This update of the variable $\nu$ at time $t_{k_i}$ is obtained by simply realizing that $\nu_i(t_{k_i})=y_{ij}(t_{k_i})-y_{ii}(t_{k_i})=y_{ii}(t_{k_i}^-)-x_i(t_{k_i}^-)=e_{ii}(t_{k_i}^-)$ (the local error just before it resets to zero), since the local model $y_{ii}$ is updated using $x_i(t_{k_i})$ and $y_{ij}(t_{k_i})=y_{ij}(t_{k_i}^-)=y_{ii}(t_{k_i}^-)$. 

Define $\nu(t_\mu)=[\nu_1^T(t_\mu)\ldots\nu_N^T(t_\mu)]^T$ and we have that
\begin{align}
  \zeta_d(t_\mu)=\nu(t_\mu) + \zeta(t_\mu). \label{eq:zetad}
\end{align}
Consider the worst case scenario (greatest difference between $\zeta_d(t_\mu)$ and $\zeta(t_\mu)$) given when all agents transmit at the same instant $t_k$ and the greatest possible delay $d \ (\geq hp_i)$ is present. Then, we have the following
\begin{align}
  \left.
	\begin{array}{l l}
 \zeta_d^T(t_k+d)\zeta_d(t_k+d)  
 \!\!\!&=\zeta^T(t_k^-)(\bar{G}^p)^T\bar{G}^p\zeta(t_k^-) \\
 &~~ + 2\zeta^T(t_k+d)\bar{G}^p\zeta(t_k^-) 
 + \zeta^T(t_k+d)\zeta(t_k+d)  \\
 &=\sum_{i=1}^N\big[e_{ii}^T(t_{k_i}^-)(G^p)^TG^pe_{ii}(t_{k_i}^-)  
  + 2e_{ii}^T(t_{k_i}+d)G^pe_{ii}(t_{k_i}^-)  \\ 
 &~~+ e_{ii}^T(t_{k_i}+d)e_{ii}(t_{k_i}+d)  \big] \\
 &\leq \sum_{i=1}^N\big[e_{ii}^T(t_{k_i}^-)(G^p)^TG^pe_{ii}(t_{k_i}^-)   
 + 2\left\|G^pe_{ii}(t_{k_i}^-)\right\|\Upsilon z_{i,M} + (\Upsilon z_{i,M})^2 \big]
	\end{array}  \right. 		\label{eq:errij-errii}
\end{align}	
where $\bar{G}=I_N\otimes G$. The local error $e_{ii}^T(t_{k_i}^-)$ represents the error just before the update instant, i.e., before it is reset to zero because of the update at time $t_{k_i}$. On the other hand, $e_{ii}(t_{k_i}+d)$ represents the error \textit{after} the update at time $t_{k_i}$ and it can only be estimated using
 \begin{align}
  \left.
	\begin{array}{l l}
 \left\|e_{ii}(t_{k_i}\!+\!d)\right\|\!\!\!\!&\leq \left\|e^{Ad}e_{ii}(t_{k_i})\right\|
   +z_{i,M}\int_0^d\left\|e^{A(d-s)}cBF\right\|ds  \\
   &\leq \Upsilon z_{i,M}
	\end{array}  \right. 		\label{eq:errtkpd}
\end{align}	
where $e_{ii}(t_{k_i})=0$ and $z_{i,M}$ represents an upper-bound on the the local control input for the time interval $t_\mu\in[t_{k_i},t_{k_i}+d]$, that is, $\left\|z_i(t_{k_i}+hp_i)\right\|\leq z_{i,M}$ for $p_i=1,...,p$.

Since the worst case is given by the maximum delay $d$ we can use \eqref{eq:errij-errii} and the current local error $e_{ii}(t_\mu)$ to bound the delayed error $e_{ij}(t_\mu+d)$ for any sampling time $t_\mu>0$, therefore the term \eqref{eq:deltadd} is used as a part of the overall threshold \eqref{eq:deltadel}. In other words, we propagate the current error $e_{ii}(t_\mu)$ forward in time using the worst case delay $d$, as if we had an event at time $t_\mu$, then we check \eqref{eq:CondDelay}. If \eqref{eq:CondDelay} holds then an event is triggered; otherwise, no event is needed and we repeat the same process at the following sampling time instant $t_{\mu+1}$.

Now, we consider the discretization error $\breve{\zeta}_d$ corresponding to the delayed state errors $\zeta_d(t_\mu)$. The discretization errors are reset at every sampling time $t_\mu$; therefore, we only need to consider the effect of these errors at time $t_k+d$ for the worst case delay $d \ (\geq hp_i)$. In other words, we need to evaluate the difference $\breve{\zeta}_d(t_k+d)=\zeta_d(t_k+d-h)-\zeta_d(t_k+d)$ as a function of the state errors $\zeta(t_k^-)$. This can be done as follows:
\begin{align}
\left.
	\begin{array}{l l}
 \breve{\zeta}_d(t_k+d)^T\breve{\zeta}_d(t_k+d)  \\
=(\zeta_d(t_k\!+\!d\!-\!h)\! - \!\zeta_d(t_k\!+\!d))^T(\zeta_d(t_k\!+\!d\!-\!h)\! - \!\zeta_d(t_k\!+\!d)) \\
=\zeta(t_k^-)^T(\bar{G}^p-\bar{G}^{p-1})^T(\bar{G}^p-\bar{G}^{p-1})\zeta(t_k^-) \\
 ~~+2(\zeta(t_k\!+\!d)-\zeta(t_k\!+\!d\!-\!h))^T(\bar{G}^p-\bar{G}^{p-1})\zeta(t_k^-) \\
 ~~+(\zeta(t_k\!+\!d)\! - \!\zeta(t_k\!+\!d\!-\!h) )^T(\zeta(t_k\!+\!d)\! - \!\zeta(t_k\!+\!d\!-\!h))  \\
=\sum_{i=1}^N\Big[e_{ii}^T(t_{k_i}^-)((G-I)G^{p-1})^T(G-I)G^{p-1}e_{ii}(t_{k_i}^-)  \\
 ~~ + 2(e_{ii}(t_{k_i}\!\!+\!d)-e_{ii}(t_{k_i}\!\!+\!d\!-\!h))^T(G-I)G^{p-1}e_{ii}(t_{k_i}^-)  \\
 ~~ + (e_{ii}(t_{k_i}\!\!+\!d)-e_{ii}(t_{k_i}\!\!+\!d\!-\!h))^T   
 (e_{ii}(t_{k_i}\!\!+\!d)-e_{ii}(t_{k_i}\!\!+\!d\!-\!h))  \big]  \\
=\sum_{i=1}^N\Big[e_{ii}^T(t_{k_i}^-)((G-I)G^{p-1})^T(G-I)G^{p-1}e_{ii}(t_{k_i}^-)  \\
 ~~ +2\left\|(G\!-\!I)G^{p-1}e_{ii}(t_{k_i}^-)\right\| (\left\|G\!-\!I\right\|\Upsilon_h \!+\! \left\|E\right\|)z_{i,M} \\
	~~ + \big((\left\|G-I\right\|\Upsilon_h +\left\|E\right\|)z_{i,M}\big)^2  \Big]
	\end{array} \label{eq:DiscDelErr} \right.
\end{align}
since we want to guarantee \eqref{eq:DiscDelErr} for any sampling time $t_\mu$ then we use the current error $e_{ii}(t_\mu)$ and we include \eqref{eq:deltacur} in the overall threshold \eqref{eq:deltadel}.

Then, the following holds:
\begin{align}
	\dot{V}&\leq x^T\bar{\mathcal{L}}x +\sum^{N}_{i=1}\left[-c_1z_i^TPBB^TPz_i + \delta_i  \right].  \label{eq:Vdotdelay2-6} 
\end{align}
Then, the local thresholds can be defined based on the local errors $e_{ii}(t_\mu)$ as in \eqref{eq:CondDelay} with $\delta_i$ given by \eqref{eq:deltadel}. When an event is triggered the error $e_{ii}$ is reset to zero and the following holds
\begin{align}
\left.
	\begin{array}{l l}
	\dot{V}	\!\!\!&\leq x^T\bar{\mathcal{L}}x +c_1\sum^{N}_{i=1}(\sigma-1)z_i^TPBB^TPz_i  
	+N\eta   \\
	&\leq x^T\bar{\mathcal{L}}x  +N\eta
	\end{array} \label{eq:Vdotdelay3} \right.
\end{align}
and the bound \eqref{eq:disagreement} on the difference between any two states follows.

Also note that 
\begin{align}
\left.
\begin{array}{l l}
  \left\|z_i(t_{k_i}+hp_i)\right\|  
	=\left\|\sum_{j\in\mathcal{N}_i}\big(x_i(t_{k_i}+hp_i) \right. 
	  \left.-x_j(t_{k_i}+hp_i)-e_{ji}(t_{k_i}+hp_i) \big)\right\|   \\
	\leq \left\|\mathcal{L}_n x(t_{k_i}+hp_i)\right\| + b_e\left\|e_{ij}(t_{k_i}+hp_i)\right\|  \\
	\leq \lambda_{\max}(\mathcal{L})\sqrt{\frac{V(t_{k_i}+hp_i)}{\lambda_{\min}(\hat{\mathcal{L}})}} 
	+ b_e\left\|z_i(t_{k_i}+hp_i)\right\| \int_0^{hp_i}\left\|e^{A(hp_i-s)}cBF\right\|ds  \\
	\leq \lambda_{\max}(\mathcal{L})\sqrt{\frac{\max_{t_\mu}\left\{V(t_\mu)\right\}}{\lambda_{\min}(\hat{\mathcal{L}})}} 
	+ b_e\left\|z_i(t_{k_i}+hp_i)\right\| \int_0^{hp_i}\left\|e^{A(hp_i-s)}cBF\right\|ds
\end{array} \right.  \label{eq:zdelb}
\end{align}
Expression \eqref{eq:zdelb} can also be written as
\begin{align}
\left.
\begin{array}{l l}
  \left\|z_i(t_{k_i}+hp_i)\right\| 
\leq \frac{\lambda_{\max}(\mathcal{L})}{1-b_e\int_0^{hp_i}\left\|e^{A(hp_i-s)}cBF\right\|ds}\sqrt{\frac{\max_{t_\mu}\left\{V(t_\mu)\right\}}{\lambda_{\min}(\hat{\mathcal{L}})}} \\
\end{array} \right.  \label{eq:zdelb2}
\end{align}
From \eqref{eq:Voft} we have that $\max_{t_\mu}\left\{V(t_\mu\right\}=V_M=\max\left\{V(0),\frac{N\eta}{\beta}\right\}$ and since the right hand side of \eqref{eq:zdelb2} increases as $p_i$ increases from $1$ to $p$, then we can conclude that
\begin{align}
  \left\|z_i(t_{k_i}+hp_i)\right\| \leq z_{i,M}= \bar{z}_i \label{eq:zdelb3}
\end{align}
for $p_i=1,...,p$.

The final task is to determine $\eta$ such that the inter-event times satisfy \eqref{eq:tkdelay} for a given $d$ that satisifies $\Upsilon<\frac{1}{b_e}$. At time $t_{k_i}$ we have that the local error is reset, i.e. $e_{ii}(t_{k_i})=0$, and no event is to be triggered by the local agent $i$ during the time interval $t_\mu\in[t_{k_i},t_{k_i}+d)$. Then, the error $e_{ii}(t_{k_i}+d)$ can be estimated as in \eqref{eq:errtkpd}. Now, the term $\delta_i(t_{k_i}+d)$ can be bounded as follows:
\begin{align}
\left.
\begin{array}{l l}
  \left\|\delta(t_{k_i}\!+\!d)\right\|  
	&\leq c_1(1+b)^2 \breve{z}_i^T(t_{k_i}\!+\!d)PBB^TP\breve{z}_i(t_{k_i}\!+\!d)     \\
	&~~+c_1(1+\frac{1}{b})\big(\left\|G^p\right\|+1\big)^2\Upsilon^2 \bar{z}_i^2  \\
	&~~+c_1(1\!+\!b)(1\!+\!\frac{1}{b})\Big(\left\|(G\!-\!I)G^{p-1}\right\|\Upsilon  
	 +\left\|(G\!-\!I)\right\|\Upsilon_h\! + \!\left\|E\right\| \Big)^2 \bar{z}_i^2 	
\end{array} \right.  \label{eq:etadel2}
\end{align}
Following similar steps as in \eqref{eq:zdelb} the term $\breve{z}_i(t_{k_i}\!+\!d)$ can be estimated as follows:
\begin{align}
\left.
\begin{array}{l l}
  \left\|\breve{z}_i(t_{k_i}\!+\!d)\right\|  
	&=\left\|\sum_{j\in\mathcal{N}_i}\big(\breve{x}_i(t_{k_i}\!+\!d) \right.  
	\left.  -\breve{x}_j(t_{k_i}\!+\!d)-\breve{e}_{ji}(t_{k_i}\!+\!d) \big)\right\|   \\
	&\leq \left\|\mathcal{L}_n \breve{x}(t_{k_i}+d)\right\| + b_e\left\|\breve{e}_{ij}(t_{k_i}+d)\right\|  \\  
	&\leq \lambda_{\max}(\mathcal{L})\sqrt{\frac{V_M}{\lambda_{\min}(\hat{\mathcal{L}})}}\big(\frac{1}{1-b_e\Upsilon}+\frac{1}{1-b_e\Upsilon_h}  \big)
\end{array} \right.  \label{eq:zdelcur}
\end{align}
Then, by the selection of $\eta$ in \eqref{eq:etadel} we can guarantee that after an event instant at $t_{k_i}$ the term $\delta_i$ does not grow enough to reach the value $\eta$ at time $t_{k_i}+d$ and threshold \eqref{eq:CondDelay} is not triggered before or at $t_{k_i}+d$. Then, the inter-event times are greater than $d$ as shown in \eqref{eq:tkdelay} and, obviously, we can guarantee that Zeno behavior does not occur at any node. $\bullet$

\begin{figure}
	\begin{center}
		\includegraphics[width=8.4cm,height=7.5cm,trim=.5cm .0cm .5cm .0cm]{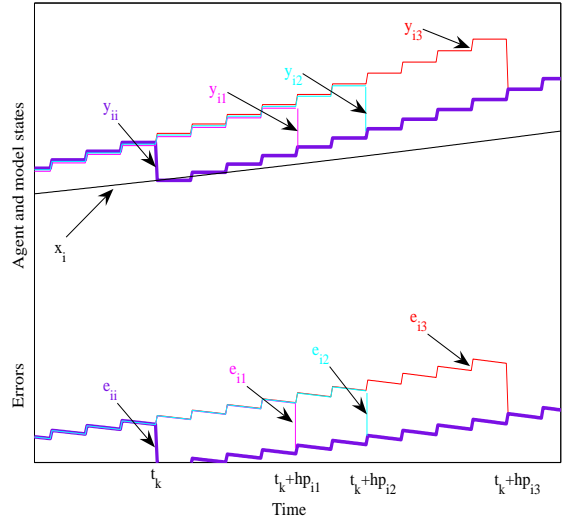}
	\caption{Relations between model states and corresponding errors when different communication delays $p_{ij}(t_{k_i})$ occur at the same update instant}
	\label{fig:DiffDel}
	\end{center}
\end{figure}

The communication delays are time-varying, that is, the value of the delays $p_i(t_{k_i})$ and $p_i(t_{k_i+1})$, for any $k_i$, are not the same in general. However, it was assumed in this section that for a given update time instant $t_{k_i}$ the neighbors of agent $i$ will all receive the measurement $x_i(t_{k_i})$ at time $t_{k_i}+hp_i(t_{k_i})$. This can be generalized to consider different delays from agent $i$ to each one of its neighbors at the same update time instant. In other words, we can consider the delays $p_{ij}(t_{k_i})$ and $p_{ij'}(t_{k_i})$ to be different in general, for $j,j' \in\mathcal{N}_i$. Thus, the models $y_{ij}$ and $y_{ij'}$ may differ during some time intervals, when one is updated sooner than the other. Fig. \ref{fig:DiffDel} shows a simple example where agents $1$, $2$, and $3$ are neighbors of agent $i$. In general, we have that $p_{i1}(t_{k_i})\neq p_{i2}(t_{k_i})\neq p_{i3}(t_{k_i})$ and the delayed models are updated at different time instants. Assuming that $p_{ij}(t_{k_i})\leq p$ for $j \in\mathcal{N}_i$, then it is clear that the same results in this section apply to the more general case under discussion since the growth on each error is still bounded by the error difference propagated using the upper-bound on the delay $d=hp\geq hp_{ij}(t_{k_i})$, for any $j \in\mathcal{N}_i$.

\section{Example} \label{sec:six}
\begin{figure}
	\begin{center}
		\includegraphics[width=8.4cm,height=7.5cm,trim=.5cm .0cm .5cm .0cm]{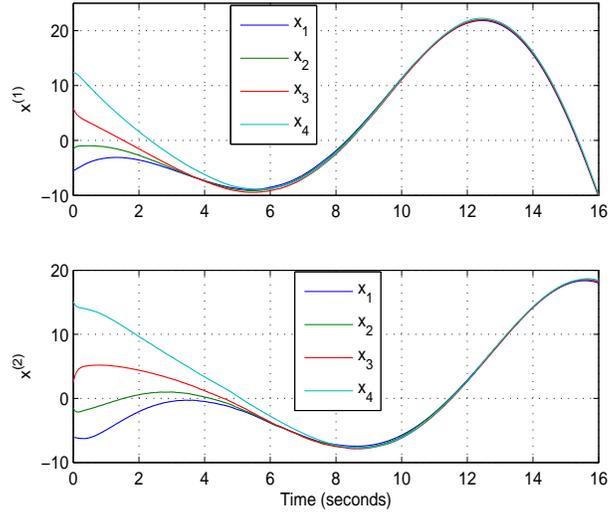}
	\caption{States of four second-order agents converging to a single trajectory}
	\label{fig:states}
	\end{center}
\end{figure}

\begin{figure}
	\begin{center}
		\includegraphics[width=8.4cm,height=7.5cm,trim=1cm .1cm .8cm .1cm]{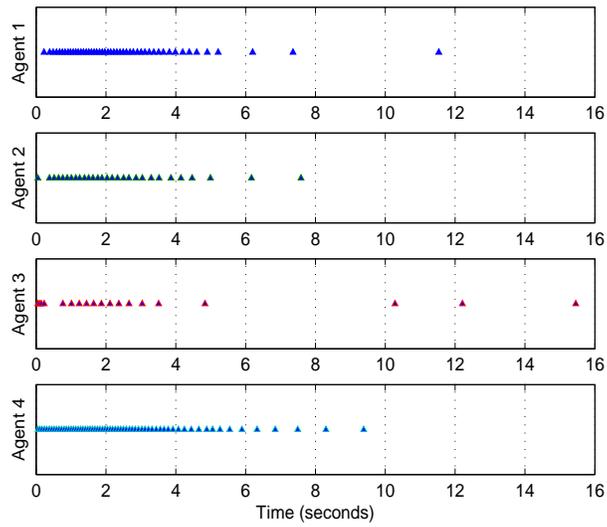}
	\caption{Broadcasting instants for each one of the four agents}
	\label{fig:comm}
	\end{center}
\end{figure}

\begin{figure}
	\begin{center}
		\includegraphics[width=8.4cm,height=7.5cm,trim=1cm .1cm .8cm .1cm]{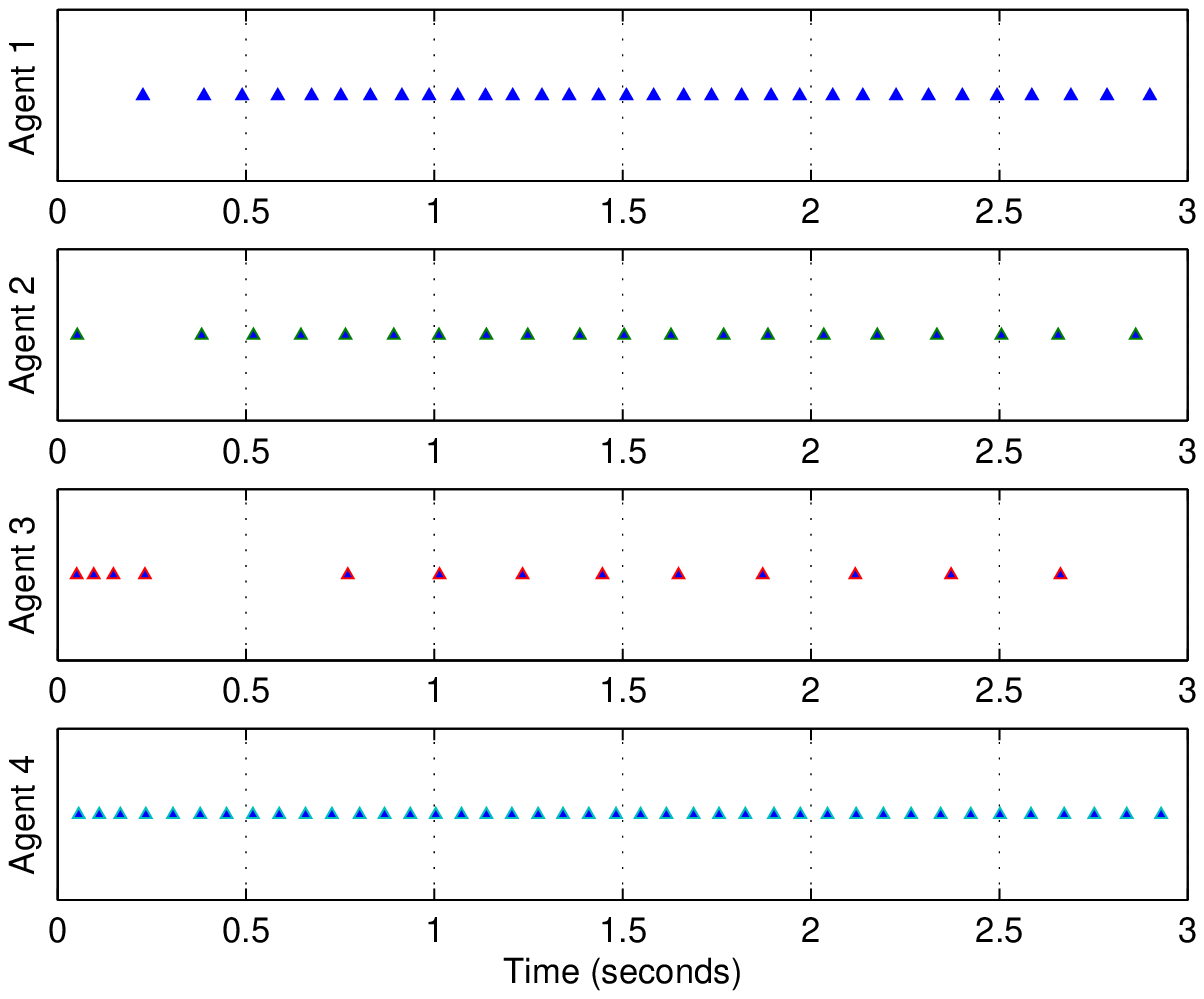}
	\caption{Broadcasting instants for each one of the four agents during the first 3 seconds of simulation}
	\label{fig:comm_zoom}
	\end{center}
\end{figure}

Consider a decentralized model-based implementation of four second-order agents ($N=4, \ n=2$) with unstable linear dynamics given by:
\[A=
\begin{bmatrix}
	0.2&-0.8\\
	0.26&0.05
\end{bmatrix}, B=
\begin{bmatrix}
	0.7\\
	-1.1
\end{bmatrix}. \]
Solving \eqref{eq:LMI} we obtain 
\[P=
\begin{bmatrix}
	0.5859&-0.1575\\
	-0.1575&0.4274
\end{bmatrix}
\]
The continuous-time state matrix is unstable with complex eigenvalues $\lambda(A)={0.125\pm .5i}$. The discretization period is $h=0.002$ seconds, the delay bound is $d=0.014$ seconds, and we compute $\eta=10.85$. The nonzero elements of the undirected adjacency matrix are $a_{12}=a_{23}=a_{34}=1$, (the corresponding symmetric elements are also equal to one). The initial condition of each agent are as follows

\[x_1(0)=
\begin{bmatrix}
	-5.5\\
	-6.1
\end{bmatrix}, x_2(0)=
\begin{bmatrix}
	-1.6\\
	-1.5
\end{bmatrix}. \] \[x_3(0)=
\begin{bmatrix}
	5.9\\
	2.5
\end{bmatrix}, x_4(0)=
\begin{bmatrix}
	12.35\\
	15.1
\end{bmatrix}. \]
Fig. \ref{fig:states} shows the response of the agents where it can be seen that the agents synchronize their states in each one of their two dimensions. Communication delays for each agent are time-varying and they take random values form the finite set: $\left\{0.010,\:0.012,\:0.014\right\}$ seconds. Fig \ref{fig:comm} shows the time instants where each agent broadcasts a measurement. Here, it can be seen that agents transmit information less frequently as they transition into a consensus state. The inter-event times for every agent are lower-bounded by the delay bound, that is, $t_{k_i+1}-t_{k_i}>0.014$ seconds. Due to a higher density of transmissions at the beginning of the simulation, a new figure, Fig. \ref{fig:comm_zoom}, has been included to show with more clarity the first $3$ seconds of the simulation and to show that the inter-event time intervals for any agent are always greater than the delay bound $d$. 

\section{Conclusions} \label{sec:seven}
Synchronization of state trajectories of linear multi-agent systems was studied in this paper. Multiple issues affecting the convergence to common trajectories were considered such as: limited sensing and actuation capabilities, limited communication, and time-varying communication delays. Event-triggered control schemes were proposed in this paper which not only provide decentralized control inputs but also allow for decentralized design of transmission instants where each agent decides, based only on local information, when to broadcast its current measurements. Then, global knowledge of communication periods and communication time instants is not needed as in sampled-data approaches. The use of discretized and decoupled models and the implementation of periodic event-triggered techniques provides a formal framework that limits actuation and sensing update rates and reduces communication. This method also provides the necessary freedom to each agent in order to determine its own broadcasting instants.

\section*{Appendix} \label{sec:Appendix}
\textit{Proof of Lemma \ref{lm:LyapFn}}. 
Sufficiency: We can express $V$ using the following:
\begin{align}
	V=\sum^N_{i=1}\sum_{j\in\mathcal{N}_i}\xi^T_iQ\left(\xi_i-\xi_j\right). \label{eq:LyapV2}
\end{align} 
Since the graph is undirected~\eqref{eq:LyapV2} can be written in the following form: 
\begin{align}
\left.
	\begin{array}{l l}
	V&=\sum^N_{i=1}\sum_{j\in\mathcal{N}_i}\frac{1}{2}\left(\xi^T_iQ\xi_i-\xi^T_iQ\xi_j\right.
	-\left.\xi^T_jQ\xi_i+\xi^T_jQ\xi_j\right)\\
  &=\sum^N_{i=1}\sum_{j\in\mathcal{N}_i}\frac{1}{2}\left(\xi_i-\xi_j\right)^TQ\left(\xi_i-\xi_j\right). 
  \end{array} \label{eq:LyapV3} \right.
\end{align} 
Since $Q>0$ and the graph is connected it is clear that if $V=0$ then $\xi_i=\xi_j$ for $i,j=1...N$. 

Necessity: Consider the following expression:
\begin{align}
   \hat{\mathcal{L}}\bar{\xi}=\left(\mathcal{L}\otimes Q\right)\bar{\xi}, \label{eq:LhatRho}
\end{align}
where $\bar{\xi}$ represents an $n-$dimensional consensus state and is given by $\bar{\xi}=\textbf{1}_N\otimes \varsigma$, where $\varsigma=\left[\varsigma_1\;\varsigma_2\;\ldots \;\varsigma_n \right]^{T}$. Then we have: 
\begin{align}
  \hat{\mathcal{L}}\bar{\xi}&=\left(\mathcal{L}\otimes Q\right)\textbf{1}_N\otimes \varsigma=\mathcal{L}\textbf{1}_N\otimes Q\varsigma=\textbf{0}_{nN}.
	\end{align}
We can conclude that if consensus is achieved then $V=\bar{\xi}^T\hat{\mathcal{L}}\bar{\xi}=0$. $\bullet$

\textit{Proof of Theorem \ref{th:consensus}}. 
Since the communication graph is undirected and connected $\mathcal{L}$ is symmetric and there exists a similarity transformation $S$ such that $\mathcal{L}_D=S^{-1}\mathcal{L}S$  is diagonal with one eigenvalue equal to zero. Define $T=S\otimes I_n$ then $\hat{\mathcal{L}}_D=T^{-1}\hat{\mathcal{L}}T=\mathcal{L}_D\otimes P$ is block diagonal with $n$ eigenvalues equal to zero.

Let us now consider the following:
\begin{align}
\left.
	\begin{array}{l l}
	T^{-1}\bar{\mathcal{L}}T&=T^{-1}\left(\hat{\mathcal{L}}A_c+A^{T}_c\hat{\mathcal{L}}\right)T\\
	&=T^{-1}\hat{\mathcal{L}}TT^{-1}A_cT+T^{-1}A^{T}_cTT^{-1}\hat{\mathcal{L}}T\\
	&=\hat{\mathcal{L}}_D\left(I_N\otimes A+c\mathcal{L}_D\otimes BF\right)
	+\left(I_N\otimes A+c\mathcal{L}_D\otimes BF\right)^{T}\hat{\mathcal{L}}_D.
	\end{array} \label{eq:LbarT} \right.
\end{align}
The term $I_N\otimes A+c\mathcal{L}_D\otimes BF$ is of the form 
	\[
\begin{bmatrix}
	A&0\\
	0&U
\end{bmatrix}
\]
In our case since $\mathcal{L}_D$ is diagonal then $U$ is block diagonal. Furthermore, each block is given by $U_i=A+c\lambda_iBF$, $i=2...N$. Then we have that~\eqref{eq:LbarT} is given by the block diagonal matrix: $diag\left\{0_n,\bar{\mathcal{L}}_2,\bar{\mathcal{L}}_3,\ldots,\bar{\mathcal{L}}_N\right\}$, where $\bar{\mathcal{L}}_i=\lambda_i\left(PA+A^TP-2c\lambda_iPBB^TP\right)$ for $i=2,...N$.
Since $c\geq1/\lambda_2$ and $P$ is the solution of~\eqref{eq:LMI} we can conclude that 
\[\lambda_i\left(PA+A^TP-2c\lambda_iPBB^TP\right)<0,\]
for $i=2...N$. We can see that $\bar{\mathcal{L}}$ has $n$ zero eigenvalues and the rest of its eigenvalues are negative. 

Consider the following:
\begin{align}
   \hat{\mathcal{L}}A_c\rho=\mathcal{L}\otimes P\left(I_N\otimes A+c\mathcal{L}\otimes BF\right)\rho, \label{eq:eigvectorsLbar}
\end{align}
where $\rho$ is an eigenvector of $\hat{\mathcal{L}}$ given by $\rho=\textbf{1}_N\otimes \varsigma$ 
associated with a zero eigenvalue of $\hat{\mathcal{L}}$ where $\varsigma=\left[\varsigma_1\;\varsigma_2\;\ldots \;\varsigma_n \right]^{T}$. Then we have:
\begin{align}
\left.
	\begin{array}{l l}
  \hat{\mathcal{L}}A_c\rho&=\left(\mathcal{L}\otimes PA+c\mathcal{L}^{2}\otimes PBF\right)\rho\\ 
	&=\mathcal{L}\textbf{1}_N\otimes PA\varsigma +c\mathcal{L}^{2}\textbf{1}_N\otimes PBF\varsigma\\
	&=\textbf{0}_{nN}.
		\end{array} \right.
\end{align}
Similarly, $A_c^T\hat{\mathcal{L}}\rho=\textbf{0}_{nN}$. Then it is clear that $\rho$ is an eigenvector of $\bar{\mathcal{L}}$ associated with a zero eigenvalue. $\bullet$


\end{document}